\documentclass{amsart}%
\usepackage{amsfonts}
\usepackage{amsmath}
\usepackage{amssymb}
\usepackage{graphicx}%
\newtheorem{theorem}{Theorem}
\theoremstyle{plain}

\newtheorem{corollary}[theorem]{Corollary}

\newtheorem{lemma}[theorem]{Lemma}

\newtheorem{proposition}[theorem]{Proposition}
\newtheorem{remark}[theorem]{Remark}
\numberwithin{equation}{section}
\numberwithin{theorem}{section}
\begin{document}
\title[Stationary Solutions]{A priori estimates of stationary solutions of an activator-inhibitor system}
\date{\today}
\author{Huiqiang Jiang}
\address{School of Mathematics, University of Minnesota, 127 Vincent Hall, 206 Church
St. S.E., Minneapolis, MN 55455}
\email{hqjiang@math.umn.edu}
\author{Wei-Ming Ni}
\address{School of Mathematics, University of Minnesota, 127 Vincent Hall, 206 Church
St. S.E., Minneapolis, MN 55455}
\email{ni@math.umn.edu}

\begin{abstract}
We consider positive solutions of the stationary Gierer-Meinhardt system
\[
\left\{
\begin{array}
[c]{lll}%
d_{1}\triangle u-u+\frac{u^{p}}{v^{q}}+\sigma=0 & \text{in} & \Omega,\\
d_{2}\triangle v-v+\frac{u^{r}}{v^{s}}=0 & \text{in} & \Omega,\\
\frac{\partial u}{\partial\nu}=\frac{\partial v}{\partial\nu}=0 & \text{on} &
\partial\Omega
\end{array}
\right.
\]
where $\triangle$ is the Laplace operator, $\Omega$ is a bounded smooth domain
in $\mathbb{R}^{n}$, $n\geq1$ and $\nu$ is the unit outer normal to
$\partial\Omega$. Under suitable conditions on the exponents $p,q,r$ and $s$,
different types of \textit{a priori} estimates are obtained, existence and
non-existence results of nontrivial solutions are derived, for both $\sigma>0$
and $\sigma=0$ cases.

\end{abstract}
\maketitle

\section{Introduction}

In 1972, following an ingenious idea of A. Turing \cite{1952Turing}, A. Gierer
and H. Meinhardt \cite{1972Gierer_Meinhardt} proposed a mathematical model for
pattern formations of spatial tissue structures of \textit{hydra} in
morphogenesis, a biological phenomenon discovered by A. Trembley in 1744
\cite{1744Trembley}. It is a system of reaction-diffusion equations of the
form
\begin{equation}
\left\{
\begin{array}
[c]{lll}%
u_{t}=d_{1}\triangle u-u+\frac{u^{p}}{v^{q}}+\sigma & \text{in} & \Omega
\times\left[  0,T\right)  ,\\
\tau v_{t}=d_{2}\triangle v-v+\frac{u^{r}}{v^{s}} & \text{in} & \Omega
\times\left[  0,T\right)  ,\\
\frac{\partial u}{\partial\nu}=\frac{\partial v}{\partial\nu}=0 & \text{on} &
\partial\Omega\times\left[  0,T\right)
\end{array}
\right. \label{GM PDE system}%
\end{equation}
where $\triangle$ is the Laplace operator, $\Omega$ is a bounded smooth domain
in $\mathbb{R}^{n}$, $n\geq1$ and $\nu$ is the unit outer normal to
$\partial\Omega$. Here $u,v$ represent respectively the concentrations of two
substances, activator and inhibitor, with diffusion rates $d_{1},d_{2}$, and
are therefore always assumed to be positive throughout this paper. The source
term $\sigma$ is a nonnegative constant representing the production of the
activator, $\tau>0$ is the response rate of $v$ to the change of $u$, and the
exponents $p,q,r,s$ are nonnegative numbers satisfying the condition
\begin{equation}
0<\frac{p-1}{r}<\frac{q}{s+1}.\label{Equation pqrs}%
\end{equation}
We remark that the response rate $\tau$ was introduced mathematically and is
an important parameter on the stability of the system.

The idea behind $\left(  \ref{GM PDE system}\right)  $ is the celebrated
\textit{diffusion-driven instability}, originally due to A. Turing
\cite{1952Turing}, which asserts that different diffusion rates could lead to
nonhomogeneous distributions of the reactants. Indeed, spike-layer stationary
solutions have been proved to exist when $\Omega$ is axially symmetric
\cite{1995Ni_Takagi} or $n=1$ \cite{1986Takagi}. When $n=2$, we refer the
readers to the discussions in \cite{2001Wei_Winter}, \cite{1999Wei_Winter} and
\cite{2002Wei_Winter}. The global existence of $\left(  \ref{GM PDE system}
\right)  $, despite partial progress made in the past $25$
years(\cite{1995Li_Chen_Chun}, \cite{1987Masuda_Takahashi}, \cite{1984Rothe}),
is only settled recently by the first author in \cite{2006Jiang} for any
positive initial data if $\frac{p-1}{r}<1$. Furthermore, if in addition
$\sigma>0$, there exists an attracting rectangle bounded away from zero and
infinity for $\left(  \ref{GM PDE system}\right)  $. On the other hand, in a
recent paper \cite{0000Ni_Suzuki_Takagi}, W. Ni, K. Suzuki and I. Takagi
completely classified the dynamics of the corresponding kinetic system
\begin{equation}
\left\{
\begin{array}
[c]{lll}%
u_{t}=-u+\frac{u^{p}}{v^{q}} & \text{in} & \left[  0,T\right)  ,\\
\tau v_{t}=-v+\frac{u^{r}}{v^{s}} & \text{in} & \left[  0,T\right)  .
\end{array}
\right. \label{GM ODE system}%
\end{equation}
In particular, it is shown that when $\frac{p-1}{r}>1$, there exist initial
values such that $u,v$ blows up in finite time. (Also see earlier results in
\cite{1995Li_Chen_Chun}.) However, the behavior of solutions to $\left(
\ref{GM PDE system}\right)  $ is not well understood in general.

Our goal here is to understand the dependence of positive steady states to
$\left(  \ref{GM PDE system}\right)  $ as the diffusion coefficients
$d_{1},d_{2}$ vary. Especially, we want to study through \textit{a priori
}estimates the existence and nonexistence of nontrivial positive stationary
patterns. When the dimension $n=1$, positive lower and upper \textit{a priori}
bounds for positive steady states of $\left(  \ref{GM PDE system}\right)  $
have been derived by I. Takagi \cite{1982Takagi}, \cite{1986Takagi} under the
general assumption $\left(  \ref{Equation pqrs}\right)  $. The method used in
\cite{1982Takagi}, \cite{1986Takagi} seems difficult to be extended to
multi-dimensional case. When $\sigma>0$ and in any space dimension, the
trivial positive lower bounds
\[
u>\sigma,\quad v>\sigma^{\frac{r}{s+1}}.
\]
immediately follow from maximum principle, and \textit{a priori} upper bounds
for H\"{o}lder norms have been obtained by W. Ni and I. Takagi
\cite{1986Ni_Takagi} using energy method under the assumption
\[
\frac{p}{q}\leq\frac{r}{s+1}\text{ and }r\geq\max\left(  p,\frac{n\left(
p-1\right)  }{2}\right)  .
\]
Also when $\sigma>0$, while studying asymptotic behavior of time-dependent
solutions to $\left(  \ref{GM PDE system}\right)  $, \textit{a priori} upper
bounds have been obtained by K. Masuda and K. Takahashi
\cite{1987Masuda_Takahashi} under the assumption
\[
\frac{p-1}{r}<\min\left\{  1,\frac{2}{n}\right\}
\]
and by the first author \cite{2006Jiang} under the assumption $\frac{p-1}
{r}<1$. When $\sigma=0$, due to the possible singularity caused by $v$ in the
denominators of the nonlinear terms, \textit{a priori} bounds usually are
harder to obtain. Nonetheless, \textit{a priori} upper bounds of H\"{o}lder
norms have been obtained in \cite{1986Ni_Takagi} under the assumption
\[
\frac{p}{q}=\frac{r}{s+1}\text{, }\frac{r}{p}>\frac{n}{2}\text{ and }
s<\frac{2}{n-2};
\]
and positive \textit{a priori} lower bounds have been obtained by M. del Pino
\cite{1994DelPino} using compactness argument under the assumption
\[
1<r<\frac{n}{n-2}\text{ and }\frac{s}{r-1}<\frac{n}{n-2}.
\]

Throughout this entire paper, we will always assume that $\left(
\ref{Equation pqrs}\right)  $ holds and use $\left(  u,v\right)  $ to denote a
smooth positive steady state of $\left(  \ref{GM PDE system}\right)  $, unless
otherwise explicitly stated, and $\left(  u^{\ast},v^{\ast}\right)  $ to
denote its unique constant steady state. We now come to our main results.

\begin{theorem}
\label{Theorem main}Suppose that $\sigma\geq0$.

\begin{enumerate}
\item[(i)] If $q<s+1$, then there exists $k_{1}>0$ depending on
$p,q,r,s,\sigma$, such that whenever $\frac{d_{2}}{d_{1}}\leq k_{1}$, we have
$u\leq u^{\ast},v\leq v^{\ast}$.

\item[(ii)] If $r<s+1$, then there exists $k_{2}>0$ depending on
$p,q,r,s,\sigma$, such that whenever $\frac{d_{2}}{d_{1}}\leq k_{2}$, we have
$u\geq u^{\ast},v\geq v^{\ast}$.

\item[(iii)] If $\max\left\{  q,r\right\}  <s+1$, then whenever $\frac{d_{2}
}{d_{1}}\leq k=\min\left\{  k_{1},k_{2}\right\}  $, we have $\left(
u,v\right)  \equiv\left(  u^{\ast},v^{\ast}\right)  $.
\end{enumerate}
\end{theorem}

\begin{remark}
The constants $k_{1},k_{2}$ and $k$ can be calculated explicitly. For example,
when $\sigma=0$ and $\left(  p,q,r,s\right)  =\left(  2,4,2,4\right)  $, the
"common source" case, we have $k_{1}=1$ and $k_{2}=11-4\sqrt{6}$, hence $k=1$.
See Theorem \ref{Theorem example} for more details.
\end{remark}

Theorem \ref{Theorem main} is new even when $n=1$. It seems interesting that
the above theorem indicates that the ratio of two diffusion rates alone can
prevent the existence of nontrivial patterns while all previously known
nonexistence results for this system require that at least one of the
diffusion rates $d_{1},d_{2}$ be suitably large. Our method also suggests that
\textit{a priori} estimates depending on $\frac{d_{2}}{d_{1}}$ are quite
natural, as the following result shows.

\begin{theorem}
\label{Theorem main 1}

\begin{enumerate}
\item[(i)] Let $\sigma=0$ and $q<s+1$. Then
\[
u\leq c\left(  1+\left(  \frac{d_{2}}{d_{1}}\right)  ^{\gamma}\right)  ,\quad
v\leq c\left(  1+\left(  \frac{d_{2}}{d_{1}}\right)  ^{\frac{r}{s+1}\gamma
}\right)
\]
where $c,\gamma$ are positive constants independent of $d_{1},d_{2}$.

\item[(ii)] Let $\sigma=0$ and
\[
\left(  \max\left\{  1,\frac{d_{2}}{d_{1}}\right\}  \right)  r<s+1.
\]
Then we have
\[
u\geq c,\quad v\geq c^{\frac{r}{s+1}},
\]
where $c\rightarrow0$ as $\frac{d_{2}}{d_{1}}\rightarrow\frac{s+1}{r}$.

\item[(iii)] Let $\sigma>0$ and $p-1<r$. Then we have
\[
u\leq c\left(  1+\left(  \frac{d_{2}}{d_{1}}\right)  ^{\gamma}\right)  ,\quad
v\leq c\left(  1+\left(  \frac{d_{2}}{d_{1}}\right)  ^{\frac{r}{s+1}\gamma
}\right)
\]
where $c,\gamma$ are positive constants independent of $d_{1},d_{2}$.
\end{enumerate}
\end{theorem}

\begin{remark}
Under the same assumption in part (iii) of Theorem \ref{Theorem main 1},
similar upper bounds have also been obtained by the first author in
\cite{2006Jiang}. However, our bounds here are more precise.
\end{remark}

A common assumption for system $\left(  \ref{GM PDE system}\right)  $ in
modeling biological pattern formation is that the activator diffuses slowly
while the inhibitor diffuses rapidly, i.e., $d_{1}$ is much smaller than
$d_{2}$. If we fix $d_{1}$ and let $d_{2}\rightarrow\infty$, formally, $v$
tends to a spatially homogeneous function $\xi=\xi\left(  t\right)  $, and
$\left(  \ref{GM PDE system}\right)  $ is reduced to the \textit{shadow
system}
\begin{equation}
\left\{
\begin{array}
[c]{lll}%
u_{t}=d_{1}\triangle u-u+\frac{u^{p}}{\xi^{q}}+\sigma & \text{in} &
\Omega\times\left[  0,T\right)  ,\\
\tau\xi_{t}=-\xi+\frac{\int_{\Omega}u^{r}\left(  x\right)  dx}{\left\vert
\Omega\right\vert \xi^{s}} & \text{in} & \Omega\times\left[  0,T\right)  ,\\
\frac{\partial u}{\partial\nu}=0 & \text{on} & \partial\Omega\times\left[
0,T\right)  .
\end{array}
\right. \label{GM shadow system}%
\end{equation}
Such formal derivation can be justified if we have \textit{a priori} estimates
independent of $d_{2}$ as $d_{2}\rightarrow\infty$.

\begin{theorem}
\label{Theorem GM estimate}

\begin{enumerate}
\item[(i)] Let $\sigma=0$ and
\[
\frac{q}{s+1}<\min\left\{  1,\frac{2}{n}\right\}  .
\]
Then we have
\[
u\leq c\left(  1+d_{1}^{-\gamma}\right)  ,\quad v\leq c\left(  1+d_{1}
^{-\frac{r}{s+1}\gamma}\right)
\]
where $c,\gamma$ are positive constants independent of $d_{1},d_{2}$.

\item[(ii)] Let $\sigma>0$, $\frac{p-1}{r}<\min\left\{  1,\frac{2}{n}\right\}
$. Then we have
\[
u\leq c\left(  1+d_{1}^{-\gamma}\right)  ,\quad v\leq c\left(  1+d_{1}
^{-\frac{r}{s+1}\gamma}\right)
\]
where $c,\gamma$ are positive constants independent of $d_{1},d_{2}$.
\end{enumerate}
\end{theorem}

The following theorem provides both lower bounds and upper bounds which are
independent of $d_{1},d_{2}$ when $d_{1},d_{2}$ are large.

\begin{theorem}
\label{Theorem GM estimate 1}Let $\sigma=0$ and $d_{1},d_{2}>\eta$ where
$\eta$ is a given positive number.

\begin{enumerate}
\item[(i)] Assume that $r<\frac{n}{n-2}$ and there exists $\delta\in\left(
0,1\right]  $ such that
\[
\frac{1-\delta}{r}+\frac{\delta}{p}<1,
\]
and
\[
\frac{\left(  \frac{1-\delta}{r}\right)  s+\frac{\delta q}{p}}{\frac
{r-1+\delta}{r}-\frac{\delta}{p}}<\frac{n}{n-2}\text{ or }\frac{\left(
\frac{1-\delta}{r}\right)  s+\frac{\delta q}{p}}{\frac{r-1+\delta}{r}
-\frac{\delta}{p}}\leq s+1.
\]
Then
\[
u\geq c_{1},\quad v\geq c_{1}^{\frac{s+1}{r}}
\]
where $c_{1}=c_{1}\left(  n,p,q,r,s,\eta\right)  $.

\item[(ii)] Assume in addition that $\frac{p-1}{r}<\min\left\{  1,\frac{2}
{n}\right\}  $. Then we have
\[
u\leq c_{2},\quad v\leq c_{2}^{\frac{r}{s+1}},
\]
where $c_{2}=c_{2}\left(  n,p,q,r,s,\eta\right)  $.
\end{enumerate}
\end{theorem}

The assumptions in part (i) of the above theorem seem complicated; however,
since $\delta$ is a free parameter, we can choose different $\delta$ to yield
a family of estimates. For example, The lower bound by M. del Pino in
\cite{1994DelPino} is contained in part (i) with $\delta=0$. Also when $n=2$,
the assumptions in part (i) hold automatically as long as we have $\left(
\ref{Equation pqrs}\right)  $.

An important consequence of the above \textit{a priori} estimates is that
$\left(  u^{\ast},v^{\ast}\right)  $ is the only steady state of $\left(
\ref{GM PDE system}\right)  $ when $d_{1}$ is suitably large. (See Theorem
\ref{Theorem nonexistence when d_1 is large}.)

\begin{theorem}
\label{Theorem nonexistence}Let $p-1<r$.

\begin{enumerate}
\item[(i)] Assume in addition $\sigma>0$. Then for any $K>0$, there exists
constant $c>0$, such that whenever $Kd_{1}\geq d_{2}$ and $d_{1}\geq c$,
$\left(  \ref{GM PDE system}\right)  $ has no nonconstant steady states.

\item[(ii)] Assume in addition $\sigma>0$ and $\frac{p-1}{r}<\frac{2}{n}$.
Then there exists constant $c>0$, such that whenever $d_{1}\geq c$, $\left(
\ref{GM PDE system}\right)  $ has no nonconstant steady states.

\item[(iii)] Assume in addition $\sigma=0$ and $n=2$. Then for any $d^{\ast
}>0$, there exists constant $c>0$, such that whenever $d_{2}\geq d^{\ast}$,
and $d_{1}\geq c$, $\left(  \ref{GM PDE system}\right)  $ has no nonconstant
steady states.
\end{enumerate}
\end{theorem}

\begin{remark}
In part (iii) of Theorem \ref{Theorem nonexistence}, assumption $n=2$ can be
replaced by the more general assumptions in Theorem
\ref{Theorem GM estimate 1}.
\end{remark}

Another application of our\textit{ a priori} estimates is the existence of
nontrivial steady states when $d_{1}$ is sufficiently small. We refer the
readers to Theorems \ref{Theorem existence} and
\ref{Theorem existence detailed} below for more details. The main idea is to
show that the Leray-Schauder degree of the associated map is nonzero in a
region excluding the trivial steady state.

Our techniques work for more general reaction-diffusion systems, but in order
to make our ideas clear, we will not pursue such generality here.

The paper is organized in the following way. We first present some basic
estimates in Section \ref{Section Preliminaries}. In Section
\ref{Section maximum principle}, we use maximum principle to establish a
priori bounds depending on $\frac{d_{2}}{d_{1}}$, especially, Theorems
\ref{Theorem main}, \ref{Theorem main 1} will be proved. In Sections
\ref{Section energy 1} and \ref{Section energy 2}, we use two different energy
methods to establish Theorems \ref{Theorem GM estimate} and
\ref{Theorem GM estimate 1}. In Section \ref{Section nonexistence}, we will
discuss nonexistence results. Finally, in Section \ref{Section existence}, we
will use topological degree theory to show the existence of nontrivial steady
states under certain situations.

\section{Preliminaries\label{Section Preliminaries}}

Let $\Omega\subset\mathbb{R}^{n}$, $n\geq1$, be a bounded smooth domain. We
consider positive stationary solutions of $\left(  \ref{GM PDE system}\right)
$
\begin{equation}
\left\{
\begin{array}
[c]{lll}%
d_{1}\triangle u-u+\frac{u^{p}}{v^{q}}+\sigma=0 & \text{in} & \Omega,\\
d_{2}\triangle v-v+\frac{u^{r}}{v^{s}}=0 & \text{in} & \Omega,\\
\frac{\partial u}{\partial\nu}=\frac{\partial v}{\partial\nu}=0 & \text{on} &
\partial\Omega,
\end{array}
\right. \label{Equation stationary}%
\end{equation}
where $d_{1},d_{2}>0$ are diffusion constants, the exponents $p,q,r,s$ are
nonnegative constants satisfying $\left(  \ref{Equation pqrs}\right)  $ and
the source term $\sigma$ is a nonnegative constant.

For any $\sigma\geq0$, $\left(  \ref{Equation stationary}\right)  $ has a
unique constant solution $\left(  u^{\ast},v^{\ast}\right)  $ such that
\begin{equation}
\left\{
\begin{array}
[c]{l}%
-u^{\ast}+\left(  u^{\ast}\right)  ^{p-\frac{qr}{s+1}}+\sigma=0,\\
v^{\ast}=\left(  u^{\ast}\right)  ^{\frac{s+1}{r}}.
\end{array}
\right. \label{Equation for u*,v*}%
\end{equation}
When $\sigma=0$, we have $\left(  u^{\ast},v^{\ast}\right)  \equiv\left(
1,1\right)  $. Furthermore, we have

\begin{proposition}
\label{Proposition for u*}For any $\sigma\geq0$, $\frac{du^{\ast}}{d\sigma}>0$
and $\lim_{\sigma\rightarrow\infty}u^{\ast}=\infty$.
\end{proposition}

\begin{proof}
Differentiating $\left(  \ref{Equation for u*,v*}\right)  _{1}$ with respect
to $\sigma$, we have
\[
-\frac{du^{\ast}}{d\sigma}+\left(  p-\frac{qr}{s+1}\right)  \left(  u^{\ast
}\right)  ^{p-1-\frac{qr}{s+1}}\frac{du^{\ast}}{d\sigma}+1=0,
\]
hence
\begin{equation}
\frac{du^{\ast}}{d\sigma}=\frac{1}{1+\left(  \frac{qr}{s+1}-p\right)  \left(
u^{\ast}\right)  ^{p-1-\frac{qr}{s+1}}}.\label{Equation
du/dsigma}%
\end{equation}
If $\frac{qr}{s+1}\geq p$, we have $\frac{du^{\ast}}{d\sigma}>0$. If
$\frac{qr}{s+1}<p$, using $\left(  \ref{Equation for u*,v*}\right)  _{1}$ to
rewrite $\left(  \ref{Equation du/dsigma}\right)  $, we have
\[
\frac{du^{\ast}}{d\sigma}=\frac{1}{\frac{qr}{s+1}-\left(  p-1\right)
+\sigma\left(  p-\frac{qr}{s+1}\right)  \left(  u^{\ast}\right)  ^{-1}}>0.
\]

\end{proof}

From now on, we assume that $\left(  u,v\right)  $ is a positive smooth
solution of $\left(  \ref{Equation stationary}\right)  $, i.e., $u,v\in
C^{\infty}\left(  \overline{\Omega}\right)  $ and $u,v>0$ in $\overline
{\Omega}$. (Actually, $u\geq0$ and $u\not \equiv 0$ implies $u,v>0$.) With
each solution $\left(  u,v\right)  $, we define the following quantities:
\[%
\begin{array}
[c]{cc}%
\bar{u}=\max_{x\in\overline{\Omega}}u, & \underline{u}=\min_{x\in
\overline{\Omega}}u,\\
\bar{v}=\max_{x\in\overline{\Omega}}v, & \underline{v}=\min_{x\in
\overline{\Omega}}v.
\end{array}
\]

First, we recall a basic convexity property of a $C^{2}$ function at its local extrema.

\begin{lemma}
\label{Lemma convexity}Let $w\in C^{2}\left(  \overline{\Omega}\right)  $
satisfy $\frac{\partial w}{\partial\nu}=0$ on $\partial\Omega$.

\begin{enumerate}
\item[(i)] If $w$ has a local maximum at $x_{1}\in\overline{\Omega}$, then
\[
\triangledown w\left(  x_{1}\right)  =0,\quad\triangle w\left(  x_{1}\right)
\leq0;
\]

\item[(ii)] If $w$ has a local minimum at $x_{2}\in\overline{\Omega}$, then
\[
\triangledown w\left(  x_{2}\right)  =0,\quad\triangle w\left(  x_{2}\right)
\geq0.
\]

\end{enumerate}
\end{lemma}

\begin{remark}
Neumann boundary condition is needed for the above lemma to hold if the local
extremum is located on $\partial\Omega$.
\end{remark}

Applying Lemma \ref{Lemma convexity} to $u,v$, we have

\begin{proposition}%
\begin{equation}
\bar{v}\leq\bar{u}^{\frac{r}{s+1}},\quad\underline{v}\geq\underline{u}
^{\frac{r}{s+1}},\quad\underline{u}\geq\frac{\underline{u}^{p}}{\bar{v}^{q}
}+\sigma,\quad\bar{u}\leq\frac{\bar{u}^{p}}{\underline{v}^{q}}+\sigma
.\label{Equation simple inequality 1}%
\end{equation}

\end{proposition}

\begin{proof}
Let $x^{\ast}\in\overline{\Omega}$ be such that $v\left(  x^{\ast}\right)
=\bar{v}$, then at $x^{\ast}$,
\[
\triangle v=\frac{1}{d_{2}}\left(  v-\frac{u^{r}}{v^{s}}\right)  \leq0,
\]
hence, $\bar{v}^{s+1}\leq u^{r}\left(  x^{\ast}\right)  \leq\bar{u}^{r}$. The
other three inequalities can be proved in the same manner.
\end{proof}

Next, we include basic energy estimates.

\begin{lemma}
\label{Lemma trivial estimates}
\[%
\begin{array}
[c]{ll}%
\int_{\Omega}\frac{u^{p}}{v^{q}}+\sigma\left\vert \Omega\right\vert
=\int_{\Omega}u, & \int_{\Omega}\frac{u^{p-1}}{v^{q}}+\sigma\int_{\Omega
}u^{-1}\leq\left\vert \Omega\right\vert ,\\
& \\
\int_{\Omega}\frac{1}{v^{q}}+\sigma\int_{\Omega}\frac{1}{u^{p}}\leq
\int_{\Omega}\frac{1}{u^{p-1}}, & \int_{\Omega}\frac{u^{r}}{v^{s}}
=\int_{\Omega}v,\\
& \\
\int_{\Omega}\frac{u^{r}}{v^{s+1}}\leq\left\vert \Omega\right\vert , &
\int_{\Omega}v^{s+1}\leq\int_{\Omega}u^{r}.
\end{array}
\]

\end{lemma}

\begin{proof}
The first identity follows from integrating $\left(  \ref{Equation stationary}
\right)  _{1}$ over $\Omega$. Next, multiplying $\left(
\ref{Equation stationary}\right)  _{1}$ with $\frac{1}{u}$ then integrating
over $\Omega$, we have
\[
\int_{\Omega}d_{1}\left\vert \triangledown u\right\vert ^{2}\frac{1}{u^{2}
}-\left\vert \Omega\right\vert +\int_{\Omega}\frac{u^{p-1}}{v^{q}}+\sigma
\int_{\Omega}u^{-1}=0,
\]
which establishes the second inequality. The other estimates can be obtained
in a similar manner.
\end{proof}

The following $L^{1}$ estimates come from standard elliptic theory.

\begin{lemma}
\label{Lemma reverse Holder}Let $\left(  u,v\right)  $ be a solution to
$\left(  \ref{Equation stationary}\right)  $. For any $0<\gamma<\frac{n}{n-2}
$, we have
\[
\int_{\Omega}u^{\gamma}\leq c\left(  1+d_{1}^{-\gamma}\right)  \left(
\int_{\Omega}u\right)  ^{\gamma},\quad\int_{\Omega}v^{\gamma}\leq c\left(
1+d_{2}^{-\gamma}\right)  \left(  \int_{\Omega}v\right)  ^{\gamma}
\]
where $c$ is a constant independent of $d_{1},d_{2}$.
\end{lemma}

\begin{proof}
Since
\[
\triangle u=\frac{1}{d_{1}}\left(  u-\frac{u^{p}}{v^{q}}-\sigma\right)  ,
\]
for any $1\leq\gamma<\frac{n}{n-2}$, we have
\begin{align*}
&  \left\Vert u\right\Vert _{\gamma}\leq c\left(  \frac{1}{d_{1}}\left\Vert
u-\frac{u^{p}}{v^{q}}-\sigma\right\Vert _{1}+\left\Vert u\right\Vert
_{1}\right) \\
\leq &  c\left(  \frac{1}{d_{1}}\left(  \left\Vert u\right\Vert _{1}
+\left\Vert \frac{u^{p}}{v^{q}}+\sigma\right\Vert _{1}\right)  +\left\Vert
u\right\Vert _{1}\right) \\
=  &  c\left(  \frac{2}{d_{1}}+1\right)  \int_{\Omega}u.
\end{align*}
And the case $0<\gamma<1$ follows from H\"{o}lder's inequality. The estimate
for $v$ can be proved in the same manner.
\end{proof}

We will need the following lemma which was proved in \cite{1994DelPino} using
Green's function approach.

\begin{lemma}
\label{Lemma del pino}\cite{1994DelPino} Let $\alpha$ be a positive constant
and $w\in C^{2}\left(  \overline{\Omega}\right)  $ be a nonnegative function
satisfying
\[
\left\{
\begin{array}
[c]{lll}%
-\triangle w+\alpha w\geq0 & \text{in} & \Omega,\\
\frac{\partial w}{\partial\nu}=0 & \text{on} & \partial\Omega.
\end{array}
\right.
\]
Then
\[
w\left(  x\right)  \geq c\int_{\Omega}w
\]
holds for any $x\in\overline{\Omega}$, here $c$ is a positive constant
depending only on $\alpha,n$ and $\Omega$.
\end{lemma}

A direct application of Lemma \ref{Lemma del pino} yields the following
estimate of solutions to $\left(  \ref{Equation stationary}\right)  $.

\begin{lemma}
\label{Lemma u,v larger than their average}
\[
u\geq c_{1}\int_{\Omega}u,\quad v\geq c_{2}\int_{\Omega}v
\]
where, for $i=1,2$, the constant $c_{i}=c_{i}\left(  n,d_{i}\right)  $ and it
can be made uniform when $d_{i}\ $is large.
\end{lemma}

\begin{proof}
Assume $d_{1}\geq\eta$, we rewrite the equation for $u$ as
\[
-\triangle u+\frac{1}{\eta}u=\frac{1}{d_{1}}\left(  \frac{u^{p}}{v^{q}}
+\sigma\right)  +\left(  \frac{1}{\eta}-\frac{1}{d_{1}}\right)  u\geq0,
\]
hence the estimate follows from Lemma \ref{Lemma del pino}.
\end{proof}

Finally, we will need the following lemma on refined Sobolev inequality.

\begin{lemma}
\label{Lemma Sobolev inequality}For any $0<\varepsilon\leq1$, $2\leq
k\leq\frac{2n}{n-2}$ if $n\geq3$ and $k\geq2$ if $n=1$ or $2$, we have for any
$w\in H^{1}\left(  \Omega\right)  $,
\[
\left(  \int_{\Omega}w^{k}dx\right)  ^{\frac{1}{k}}\leq C\varepsilon^{\frac
{n}{k}-\frac{n}{2}}\left(  \int_{\Omega}\left(  \varepsilon^{2}\left\vert
\triangledown w\right\vert ^{2}+\left\vert w\right\vert ^{2}\right)
dx\right)  ^{\frac{1}{2}},
\]
where $C$ is a positive constant independent of $\varepsilon$.
\end{lemma}

\begin{proof}
Let $w_{\varepsilon}\left(  y\right)  =w\left(  \varepsilon y\right)  $, and
$x=\varepsilon y$, then
\begin{align*}
&  \int_{\Omega}\left(  \varepsilon^{2}\left\vert \triangledown w\right\vert
^{2}+\left\vert w\right\vert ^{2}\right)  dx=\varepsilon^{n}\int_{\frac
{1}{\varepsilon}\Omega}\left(  \left\vert \triangledown w_{\varepsilon
}\right\vert ^{2}+\left\vert w_{\varepsilon}\right\vert ^{2}\right)  dy\\
\geq &  C\varepsilon^{n}\left(  \int_{\frac{1}{\varepsilon}\Omega
}w_{\varepsilon}^{k}dy\right)  ^{\frac{2}{k}}=C\varepsilon^{n-\frac{2n}{k}
}\left(  \int_{\Omega}w^{k}dx\right)  ^{\frac{2}{k}},
\end{align*}
here $C$ depends on $n,k,$ and cone property of $\frac{1}{\varepsilon}\Omega$.
Hence, for any $\varepsilon\leq1$, we have
\[
\int_{\Omega}\left(  \varepsilon^{2}\left\vert \triangledown w\right\vert
^{2}+\left\vert w\right\vert ^{2}\right)  dx\geq C\varepsilon^{n-\frac{2n}{k}
}\left(  \int_{\Omega}w^{k}dx\right)  ^{\frac{2}{k}}
\]
where $C$ is independent of $\varepsilon$.
\end{proof}

\section{Maximum Principle\label{Section maximum principle}}

In this section, we will deduce \textit{a priori} bounds for positive
solutions of $\left(  \ref{Equation stationary}\right)  $ which depend on
$\frac{d_{2}}{d_{1}}$.

First, we apply Lemma \ref{Lemma convexity} to $\frac{u}{v^{\lambda}}$, where
$\lambda$ is any real number.

\begin{lemma}
\label{Lemma local maximum}Let $0\leq\lambda\leq1$. Then
\begin{equation}
1-\frac{u^{p-1}}{v^{q}}-\sigma u^{-1}-\frac{\lambda d_{1}}{d_{2}}\left(
1-\frac{u^{r}}{v^{s+1}}\right)  \leq
0\label{Equation inequality at local maximum}%
\end{equation}
holds at any point $x^{\ast}\in\overline{\Omega}$ where $\frac{u}{v^{\lambda}
}$ achieves its local maximum.
\end{lemma}

\begin{proof}
If $\frac{u}{v^{\lambda}}$ has a local maximum at $x^{\ast}$, then at
$x^{\ast}$
\begin{equation}
\triangledown\frac{u}{v^{\lambda}}=\frac{\triangledown u}{v^{\lambda}}
-\lambda\frac{u}{v^{\lambda+1}}\triangledown
v=0,\label{Equation gradient of u/v^lambda 0}%
\end{equation}
and
\begin{equation}
\triangle\frac{u}{v^{\lambda}}=\frac{\triangle u}{v^{\lambda}}-\lambda\frac
{u}{v^{\lambda+1}}\triangle v-2\lambda\frac{\triangledown u\cdot\triangledown
v}{v^{\lambda+1}}+\lambda\left(  \lambda+1\right)  \frac{u}{v^{\lambda+2}
}\left\vert \triangledown v\right\vert ^{2}\leq
0.\label{Equation laplacian of u/v^lambda}%
\end{equation}
From $\left(  \ref{Equation gradient of u/v^lambda 0}\right)  $, we have
\begin{equation}
\triangledown u=\lambda\frac{u}{v}\triangledown v.\label{Equation
gradient of u/v^lambda}%
\end{equation}
Hence dividing $\left(  \ref{Equation laplacian of u/v^lambda}\right)  $ by
$\frac{v^{\lambda}}{u}$ and using $\left(  \ref{Equation stationary}\right)  $
together with $\left(  \ref{Equation gradient of u/v^lambda}\right)  $, we
deduce
\[
\frac{1}{d_{1}}\left(  1-\frac{u^{p-1}}{v^{q}}-\sigma u^{-1}\right)
-\frac{\lambda}{d_{2}}\left(  1-\frac{u^{r}}{v^{s+1}}\right)  -\lambda\left(
\lambda-1\right)  \frac{\left\vert \triangledown v\right\vert ^{2}}{v^{2}}
\leq0.
\]
Since $\lambda\left(  \lambda-1\right)  \leq0$, $\left(
\ref{Equation inequality at local maximum}\right)  $ follows.
\end{proof}

Similarly, we have

\begin{lemma}
\label{Lemma local minimum}Let $\lambda\leq0$ or $\lambda\geq1$. Then
\begin{equation}
1-\frac{u^{p-1}}{v^{q}}-\sigma u^{-1}-\frac{\lambda d_{1}}{d_{2}}\left(
1-\frac{u^{r}}{v^{s+1}}\right)  \geq0
\end{equation}
holds at any point $x^{\ast}\in\overline{\Omega}$ where $\frac{u}{v^{\lambda}
}$ achieves its local minimum.
\end{lemma}

Our strategy is to use the bounds of $\frac{u}{v^{\lambda}}$ to control $u$
and $v$ because we have

\begin{lemma}
\label{Lemma from u/v^lambda to u,v}Let $0<\lambda<\frac{s+1}{r}$, then
\[
\inf_{\Omega}\frac{u}{v^{\lambda}}\leq\underline{u}^{1-\frac{r}{s+1}\lambda
}\leq\underline{v}^{\frac{s+1}{r}-\lambda}\leq\bar{v}^{\frac{s+1}{r}-\lambda
}\leq\bar{u}^{1-\frac{r}{s+1}\lambda}\leq\sup_{\Omega}\frac{u}{v^{\lambda}}.
\]

\end{lemma}

\begin{proof}
Let $x\in\overline{\Omega}$ be such that $u\left(  x\right)  =\bar{u}$, then
we have
\[
\bar{u}=u\left(  x\right)  =\frac{u\left(  x\right)  }{v^{\lambda}\left(
x\right)  }v^{\lambda}\left(  x\right)  \leq\left(  \sup_{\Omega}\frac
{u}{v^{\lambda}}\right)  \bar{v}^{\lambda}\leq\left(  \sup_{\Omega}\frac
{u}{v^{\lambda}}\right)  \bar{u}^{\frac{r}{s+1}\lambda},
\]
hence
\[
\bar{u}^{1-\frac{r}{s+1}\lambda}\leq\sup_{\Omega}\frac{u}{v^{\lambda}}.
\]
Similarly, we can prove
\[
\underline{u}^{1-\frac{r}{s+1}\lambda}\geq\inf_{\Omega}\frac{u}{v^{\lambda}}.
\]
The inequalities in between follow from $\left(
\ref{Equation simple inequality 1}\right)  $.
\end{proof}

Next, we introduce a family of functions
\begin{equation}
f_{\sigma}\left(  \lambda\right)  =\frac{\lambda\left(  s+1-\lambda r\right)
}{\left(  q-\lambda\left(  p-1\right)  \right)  \left(  u^{\ast}\right)
^{-\frac{qr-\left(  p-1\right)  \left(  s+1\right)  }{s+1}}+\lambda
\sigma\left(  u^{\ast}\right)  ^{-1}}\label{Equation f_sigma}%
\end{equation}
which are continuous and nonnegative in $\left[  0,\frac{s+1}{r}\right]  $. It
is easy to check that
\[
f_{\sigma}\left(  0\right)  =f_{\sigma}\left(  \frac{s+1}{r}\right)  =0
\]
and $f_{\sigma}$ has a unique critical point in $\left[  0,\frac{s+1}
{r}\right]  $ which is a local maximum.

Using $\left(  \ref{Equation for u*,v*}\right)  _{1}$, we have
\[
f_{\sigma}\left(  \lambda\right)  =\frac{\lambda\left(  s+1-\lambda r\right)
}{\left(  q-\lambda p\right)  \left(  u^{\ast}\right)  ^{-\frac{qr-\left(
p-1\right)  \left(  s+1\right)  }{s+1}}+\lambda}.
\]
Proposition \ref{Proposition for u*} implies that $f_{\sigma}\left(
\lambda\right)  $ is monotone in $\sigma$ for any given $\lambda\in\left(
0,\frac{s+1}{r}\right)  $, and
\[
\lim_{\sigma\rightarrow\infty}f_{\sigma}\left(  \lambda\right)  =s+1-\lambda
r.
\]

When $\sigma=0$, we have $u^{\ast}=1$, and the function $f_{\sigma}$ has a
simple form
\[
f_{0}\left(  \lambda\right)  =\frac{\lambda\left(  s+1-\lambda r\right)
}{q-\lambda\left(  p-1\right)  }.
\]
The unique critical point of $f_{0}$ in $\left(  0,\frac{s+1}{r}\right)  $ can
be calculated
\[
\lambda^{\ast}=\frac{q}{p-1}\left(  1-\sqrt{1-\frac{\left(  p-1\right)
\left(  s+1\right)  }{qr}}\right)  =\frac{\frac{s+1}{r}}{1+\sqrt
{1-\frac{\left(  p-1\right)  \left(  s+1\right)  }{qr}}}.
\]

We define the quantities
\[%
\begin{array}
[c]{ll}%
a=\frac{q-\lambda\left(  p-1\right)  }{s+1-\lambda r}, & b=\frac{s+1-\lambda
r}{q-\lambda\left(  p-1\right)  }=\frac{1}{a},\\
& \\
a_{0}=\frac{\lambda}{s+1-\lambda r}, & b_{0}=\frac{s+1-\lambda r}{\lambda
}=\frac{1}{a_{0}}%
\end{array}
\]
which will appear frequently in our proofs. The quantities $\frac{u^{p-1}
}{v^{q}}$ and $\frac{u^{r}}{v^{s+1}}$ are related through $\frac{u}
{v^{\lambda}}$ in the following identities
\[
\frac{u^{p-1}}{v^{q}}=\left(  \frac{u^{r}}{v^{s+1}}\right)  ^{a}\left(
\left(  \frac{u}{v^{\lambda}}\right)  ^{-\frac{qr-\left(  p-1\right)  \left(
s+1\right)  }{\left(  s+1\right)  -q-\lambda\left(  r-\left(  p-1\right)
\right)  }}\right)  ^{1-a},
\]
\[
\frac{u^{r}}{v^{s+1}}=\left(  \frac{u^{p-1}}{v^{q}}\right)  ^{b}\left(
\left(  \frac{u}{v^{\lambda}}\right)  ^{-\frac{qr-\left(  p-1\right)  \left(
s+1\right)  }{\left(  s+1\right)  -q-\lambda\left(  r-\left(  p-1\right)
\right)  }}\right)  ^{1-b}.
\]
When $\sigma>0$, we need the following identities to connect $\frac{u^{r}
}{v^{s+1}}$ and $\frac{1}{u}$,
\[
\frac{1}{u}=\left(  \frac{u^{r}}{v^{s+1}}\right)  ^{a_{0}}\left(  \left(
\frac{u}{v^{\lambda}}\right)  ^{-\frac{qr-\left(  p-1\right)  \left(
s+1\right)  }{\left(  s+1\right)  -q-\lambda\left(  r-\left(  p-1\right)
\right)  }}\right)  ^{1-a_{0}},
\]
\[
\frac{u^{r}}{v^{s+1}}=\left(  \frac{1}{u}\right)  ^{b_{0}}\left(  \left(
\frac{u}{v^{\lambda}}\right)  ^{-\frac{qr-\left(  p-1\right)  \left(
s+1\right)  }{\left(  s+1\right)  -q-\lambda\left(  r-\left(  p-1\right)
\right)  }}\right)  ^{1-b_{0}}.
\]
When $0<\lambda<\frac{s+1}{r}$, using quantities $a$ and $a_{0}$, we have
\[
f_{\sigma}\left(  \lambda\right)  =\frac{\lambda}{a\left(  u^{\ast}\right)
^{-\frac{qr-\left(  p-1\right)  \left(  s+1\right)  }{s+1}}+a_{0}\sigma\left(
u^{\ast}\right)  ^{-1}}=\frac{\lambda}{\left(  a-a_{0}\right)  \left(
u^{\ast}\right)  ^{-\frac{qr-\left(  p-1\right)  \left(  s+1\right)  }{s+1}
}+a_{0}}.
\]
And when $\sigma=0$, we have
\[
f_{0}\left(  \lambda\right)  =\frac{\lambda}{a}.
\]

Our first result is the upper bounds when $\sigma=0$.

\begin{theorem}
\label{Theorem upper bounds}Assume $\sigma=0$ and $q<s+1$. Then there exists
positive constants $c\ $and $\gamma$ independent of $d_{1},d_{2}$ such that
\[
u\leq c\left(  1+\left(  \frac{d_{2}}{d_{1}}\right)  ^{\gamma}\right)  \text{
and }v\leq c\left(  1+\left(  \frac{d_{2}}{d_{1}}\right)  ^{\frac{r}
{s+1}\gamma}\right)  .
\]
Furthermore, if $\frac{d_{2}}{d_{1}}\in f_{0}\left(  \Lambda_{1}\right)  $,
where
\[
\Lambda_{1}=\left(  0,1\right]  \cap\left(  0,\frac{s+1-q}{r-\left(
p-1\right)  }\right]  ,
\]
then $u\leq1,v\leq1$.
\end{theorem}

\begin{proof}
Let $\lambda\in\tilde{\Lambda}_{1}\cap\left(  0,\frac{d_{2}}{d_{1}}\right)  $,
where
\[
\tilde{\Lambda}_{1}=\left(  0,1\right]  \cap\left(  0,\frac{s+1-q}{r-\left(
p-1\right)  }\right)  .
\]
Since $0<\lambda\leq1$, we have from Lemma \ref{Lemma local maximum}, at any
point $x^{\ast}\in\overline{\Omega}$ where $\frac{u}{v^{\lambda}}$ achieves
its maximum,
\begin{equation}
1-\frac{\lambda d_{1}}{d_{2}}\leq\frac{u^{p-1}}{v^{q}}-\frac{\lambda d_{1}
}{d_{2}}\frac{u^{r}}{v^{s+1}}.\label{Equation inequality at
local maximum 1}%
\end{equation}
Since
\[
0<\frac{p-1}{r}<\frac{q}{s+1}<1,
\]
we have
\[
\lambda<\frac{s+1-q}{r-\left(  p-1\right)  }<\frac{s+1}{r}<\frac{q}{p-1},
\]
hence
\[
a=\frac{q-\lambda\left(  p-1\right)  }{s+1-\lambda r}\in\left(  0,1\right)  .
\]
Applying Young's inequality, we have
\begin{align*}
&  \frac{u^{p-1}}{v^{q}}=\left(  \frac{u^{r}}{v^{s+1}}\right)  ^{a}\left(
\left(  \frac{u}{v^{\lambda}}\right)  ^{-\frac{qr-\left(  p-1\right)  \left(
s+1\right)  }{s+1-\lambda r-\left(  q-\lambda\left(  p-1\right)  \right)  }
}\right)  ^{1-a}\\
=  &  \left(  \frac{1}{a}\frac{\lambda d_{1}}{d_{2}}\frac{u^{r}}{v^{s+1}
}\right)  ^{a}\left(  \left[  \frac{1}{a}\frac{\lambda d_{1}}{d_{2}}\right]
^{-\frac{a}{1-a}}\left(  \frac{u}{v^{\lambda}}\right)  ^{-\frac{qr-\left(
p-1\right)  \left(  s+1\right)  }{s+1-\lambda r-\left(  q-\lambda\left(
p-1\right)  \right)  }}\right)  ^{1-a}\\
\leq &  \frac{\lambda d_{1}}{d_{2}}\frac{u^{r}}{v^{s+1}}+\frac{1-a}{\left[
\frac{1}{a}\frac{\lambda d_{1}}{d_{2}}\right]  ^{\frac{a}{1-a}}}\left(
\frac{u}{v^{\lambda}}\right)  ^{-\frac{qr-\left(  p-1\right)  \left(
s+1\right)  }{s+1-\lambda r-\left(  q-\lambda\left(  p-1\right)  \right)  }}\\
\leq &  \frac{\lambda d_{1}}{d_{2}}\frac{u^{r}}{v^{s+1}}+\frac{1-a}{\left[
\frac{1}{a}\frac{\lambda d_{1}}{d_{2}}\right]  ^{\frac{a}{1-a}}}\left(
\bar{u}\right)  ^{-\frac{qr-\left(  p-1\right)  \left(  s+1\right)  }{\left(
s+1\right)  \left(  1-a\right)  }}%
\end{align*}
where we used Lemma \ref{Lemma from u/v^lambda to u,v} in the last inequality.
Combining this with $\left(  \ref{Equation inequality at local maximum 1}
\right)  $, we have
\[
1-\frac{\lambda d_{1}}{d_{2}}\leq\frac{1-a}{\left[  \frac{1}{a}\frac{\lambda
d_{1}}{d_{2}}\right]  ^{\frac{a}{1-a}}}\left(  \bar{u}\right)  ^{-\frac
{qr-\left(  p-1\right)  \left(  s+1\right)  }{\left(  s+1\right)  \left(
1-a\right)  }}.
\]
Since
\[
1-\frac{\lambda d_{1}}{d_{2}}>0,
\]
we deduce
\begin{equation}
\bar{u}\leq\left[  \frac{\left(  1-a\right)  ^{1-a}a^{a}}{\left(
1-\frac{\lambda d_{1}}{d_{2}}\right)  ^{1-a}\left(  \frac{\lambda d_{1}}
{d_{2}}\right)  ^{a}}\right]  ^{\frac{s+1}{qr-\left(  p-1\right)  \left(
s+1\right)  }}.\label{Equation upper bounds for
u when sigma=0}%
\end{equation}

Next we prove the optimal bounds, if
\[
\frac{d_{2}}{d_{1}}\in f_{0}\left(  \Lambda_{1}\right)  ,
\]
then we have for some $\lambda\in\Lambda_{1}$,
\[
\frac{d_{2}}{d_{1}}=f_{0}\left(  \lambda\right)  =\frac{\lambda}{a}.
\]
We first assume $\lambda\in\tilde{\Lambda}_{1}$, then we have $a\in\left(
0,1\right)  $ and
\[
\lambda=\frac{d_{2}}{d_{1}}\cdot a<\frac{d_{2}}{d_{1}},
\]
hence $\left(  \ref{Equation upper bounds for u when sigma=0}\right)  $ holds
and the bound becomes $\bar{u}\leq1$. If $\lambda\not \in \tilde{\Lambda}_{1}
$, then $\lambda=\frac{s+1-q}{r-\left(  p-1\right)  }$ and $a=1$, hence
\[
\frac{\lambda d_{1}}{d_{2}}=1.
\]
Furthermore, there exists a sequence $\left\{  \lambda_{k}\right\}
\subset\tilde{\Lambda}_{1}$ such that $\lambda<\lambda_{k}$ and
\[
\lambda=\lim_{k\rightarrow\infty}\lambda_{k}.
\]
Let
\[
a_{k}=\frac{q-\lambda_{k}\left(  p-1\right)  }{s+1-\lambda_{k}r},
\]
then $a_{k}\in\left(  0,1\right)  $ and
\[
\lim_{k\rightarrow\infty}a_{k}=1.
\]
For each $k$, we have
\[
\bar{u}\leq\left[  \frac{\left(  1-a_{k}\right)  ^{1-a_{k}}a_{k}^{a_{k}}
}{\left(  1-\frac{\lambda_{k}d_{1}}{d_{2}}\right)  ^{1-a_{k}}\left(
\frac{\lambda_{k}d_{1}}{d_{2}}\right)  ^{a_{k}}}\right]  ^{\frac
{s+1}{qr-\left(  p-1\right)  \left(  s+1\right)  }}.
\]
Now
\begin{align*}
&  \lim_{k\rightarrow\infty}\frac{\left(  1-a_{k}\right)  ^{1-a_{k}}
a_{k}^{a_{k}}}{\left(  1-\frac{\lambda_{k}d_{1}}{d_{2}}\right)  ^{1-a_{k}
}\left(  \frac{\lambda_{k}d_{1}}{d_{2}}\right)  ^{a_{k}}}\\
=  &  \lim_{k\rightarrow\infty}\left(  \frac{1-a_{k}}{1-\frac{\lambda_{k}
d_{1}}{d_{2}}}\right)  ^{1-a_{k}}=\lim_{k\rightarrow\infty}\left(
\frac{s+1-q}{s+1-\lambda_{k}r}\right)  ^{1-a_{k}}=1,
\end{align*}
hence $\bar{u}\leq1$. And the optimal bound $\bar{v}\leq1$ follows from
$\left(  \ref{Equation simple inequality 1}\right)  $. The set $\Lambda_{1}$
is a nonempty interval with left end point zero, from the property of function
$f_{0}$, the set $f_{0}\left(  \Lambda_{1}\right)  $ is a nonempty interval of
the form $\left(  0,k_{1}\right]  $ for some constant $k_{1}$ depending on
$p,q,r,s$. Inequality
\[
u\leq c\left(  1+\left(  \frac{d_{2}}{d_{1}}\right)  ^{\gamma}\right)
\]
automatically holds if $\frac{d_{2}}{d_{1}}\leq k_{1}$.$\ $When $\frac{d_{2}
}{d_{1}}\geq k_{1}$, it can be deduced from $\left(
\ref{Equation upper bounds for u when sigma=0}\right)  $ by taking
\[
\lambda=\frac{1}{2}\min\left\{  1,k_{1},\frac{s+1-q}{r-\left(  p-1\right)
}\right\}  \in\tilde{\Lambda}_{1}\cap\left(  0,\frac{d_{2}}{d_{1}}\right)  .
\]
Finally
\[
v\leq c\left(  1+\left(  \frac{d_{2}}{d_{1}}\right)  ^{\frac{r}{s+1}\gamma
}\right)
\]
follows from $\left(  \ref{Equation simple inequality 1}\right)  $.
\end{proof}

\begin{remark}
It is easy to see that for any $\tau,t\in\left(  0,1\right)  $, we always
have
\[
\tau^{t}\left(  1-\tau\right)  ^{1-t}\leq t^{t}\left(  1-t\right)  ^{1-t},
\]
and the equality holds if and only if $\tau=t$.
\end{remark}

\begin{remark}
$f_{0}\left(  \Lambda_{1}\right)  $ is a nonempty interval of the form
$\left(  0,k_{1}\right]  $ for some constant $k_{1}$ depending only on
$p,q,r,s$. If $\frac{s+1-q}{r-\left(  p-1\right)  }>1$, then we have
$\Lambda_{1}=\left(  0,1\right]  $ and
\[
f_{0}\left(  \Lambda_{1}\right)  \supset\left(  0,\frac{s+1-r}{q-\left(
p-1\right)  }\right]  .
\]
And if $\frac{s+1-q}{r-\left(  p-1\right)  }\leq1$, then we have
\[
\Lambda_{1}=\left(  0,\frac{s+1-q}{r-\left(  p-1\right)  }\right]
\]
and
\[
f_{0}\left(  \Lambda_{1}\right)  \supset\left(  0,\frac{s+1-q}{r-\left(
p-1\right)  }\right]  .
\]
Especially, if $\left(  p,q,r,s\right)  =\left(  2,4,2,4\right)  $, we have
$\Lambda_{1}=\left(  0,1\right]  $, and since the critical point of $f_{0}$ in
$\left(  0,\frac{s+1}{r}\right)  $,
\[
\lambda^{\ast}=\frac{\frac{s+1}{r}}{1+\sqrt{1-\frac{\left(  p-1\right)
\left(  s+1\right)  }{qr}}}=\frac{\frac{5}{2}}{1+\sqrt{1-\frac{5}{8}}}
=\frac{10}{4+\sqrt{6}}>1,
\]
we have $f_{0}\left(  \Lambda_{1}\right)  =\left(  0,1\right]  $.
\end{remark}

The optimal bounds $u\leq1,v\leq1$ when $\frac{d_{2}}{d_{1}}$ is sufficiently
small indicate that \textit{a priori} estimates depending on $d_{1},d_{2}$ in
terms of $\frac{d_{2}}{d_{1}}$ could be natural. Such estimates are new even
in one dimensional case.

Next, we consider lower bounds of stationary solutions when $\sigma=0$.

\begin{theorem}
\label{Theorem lower bounds}Assume $\sigma=0$ and
\begin{equation}
\left(  \max\left\{  1,\frac{d_{2}}{d_{1}}\right\}  \right)
r<s+1.\label{Equation assumption lower}%
\end{equation}
Then
\[
u\geq c,\quad v\geq c^{\frac{r}{s+1}}
\]
where $c$ is a positive constant depending on $p,q,r,s$ and $\frac{d_{2}
}{d_{1}}$ and satisfies
\[
\lim_{\frac{d_{2}}{d_{1}}\rightarrow\frac{s+1}{r}}c=0.
\]
Furthermore, if $\frac{d_{2}}{d_{1}}\in f_{0}\left(  \Lambda_{2}\right)  $,
where
\[
\Lambda_{2}=\left\{  \lambda\in\left[  1,\frac{s+1}{r}\right)  :\frac
{s+1-\lambda r}{q-\lambda\left(  p-1\right)  }\leq1\right\}  .
\]
then $u\geq1,v\geq1$.
\end{theorem}

\begin{proof}
Let $\lambda\in\left(  \frac{d_{2}}{d_{1}},\infty\right)  \cap\tilde{\Lambda
}_{2}$, where
\[
\tilde{\Lambda}_{2}=\left\{  \lambda\in\left[  1,\frac{s+1}{r}\right)
:\frac{s+1-\lambda r}{q-\lambda\left(  p-1\right)  }<1\right\}  .
\]
Since $\lambda\geq1$, we have from Lemma \ref{Lemma local minimum},
\begin{equation}
\frac{\lambda d_{1}}{d_{2}}-1\leq\frac{\lambda d_{1}}{d_{2}}\frac{u^{r}
}{v^{s+1}}-\frac{u^{p-1}}{v^{q}}\label{Equation inequality at local minimum 1}%
\end{equation}
holds at any point $x^{\ast}\in\overline{\Omega}$ where $\frac{u}{v^{\lambda}
}$ achieves its minimum. Since
\[
\lambda<\frac{s+1}{r}<\frac{q}{p-1},
\]
we have
\[
b=\frac{s+1-\lambda r}{q-\lambda\left(  p-1\right)  }\in\left(  0,1\right)  .
\]
Applying Young's inequality, we have
\begin{align*}
&  \frac{\lambda d_{1}}{d_{2}}\frac{u^{r}}{v^{s+1}}\\
=  &  \left[  \frac{1}{b}\frac{u^{p-1}}{v^{q}}\right]  ^{b}\left(  b^{\frac
{b}{1-b}}\left(  \frac{\lambda d_{1}}{d_{2}}\right)  ^{\frac{1}{1-b}}\left(
\frac{u}{v^{\lambda}}\right)  ^{\frac{qr-\left(  p-1\right)  \left(
s+1\right)  }{q-\lambda\left(  p-1\right)  -\left(  s+1-\lambda r\right)  }
}\right)  ^{1-b}\\
\leq &  \frac{u^{p-1}}{v^{q}}+\left(  1-b\right)  b^{\frac{b}{1-b}}\left(
\frac{\lambda d_{1}}{d_{2}}\right)  ^{\frac{1}{1-b}}\left(  \frac
{u}{v^{\lambda}}\right)  ^{\frac{qr-\left(  p-1\right)  \left(  s+1\right)
}{q-\lambda\left(  p-1\right)  -\left(  s+1-\lambda r\right)  }}\\
\leq &  \frac{u^{p-1}}{v^{q}}+\left(  1-b\right)  b^{\frac{b}{1-b}}\left(
\frac{\lambda d_{1}}{d_{2}}\right)  ^{\frac{1}{1-b}}\underline{u}
^{\frac{qr-\left(  p-1\right)  \left(  s+1\right)  }{s+1}\cdot\frac{b}{1-b}},
\end{align*}
where we have used Lemma \ref{Lemma from u/v^lambda to u,v} in the last
inequality. Combining $\left(  \ref{Equation inequality at local maximum 1}
\right)  $, we have
\begin{equation}
\underline{u}\geq\left[  \frac{\left(  \frac{d_{2}}{\lambda d_{1}}\right)
^{b}\left(  1-\frac{d_{2}}{\lambda d_{1}}\right)  ^{1-b}}{b^{b}\left(
1-b\right)  ^{1-b}}\right]  ^{\frac{1}{b}\cdot\frac{s+1}{qr-\left(
p-1\right)  \left(  s+1\right)  }}%
\label{Equation lower bounds for u when sigma=0}%
\end{equation}
which yields a lower bound for $u$. Lower bound for $v$ follows from $\left(
\ref{Equation simple inequality 1}\right)  $.

Next, if $\frac{d_{2}}{d_{1}}\in f_{0}\left(  \Lambda_{2}\right)  $, then we
have for some $\lambda\in\Lambda_{2}$,
\[
\frac{d_{2}}{d_{1}}=f_{0}\left(  \lambda\right)  =\frac{\lambda\left(
s+1-\lambda r\right)  }{q-\lambda\left(  p-1\right)  }=\lambda b.
\]
If $\lambda\in\tilde{\Lambda}_{2}$, then $b\in\left(  0,1\right)  $ and
$\lambda>\frac{d_{2}}{d_{1}}$, hence $\left(
\ref{Equation lower bounds for u when sigma=0}\right)  $ holds and becomes
$\underline{u}\geq1$. If $\lambda\not \in \tilde{\Lambda}_{2}$, then
\[
b=\frac{s+1-\lambda r}{q-\lambda\left(  p-1\right)  }=1
\]
and there exist $\lambda_{k}\in\tilde{\Lambda}_{2}$ such that $\lambda
_{k}>\lambda$ and
\[
\lim_{k\rightarrow\infty}\lambda_{k}=\lambda.
\]
Now for each $k$, we have
\[
\underline{u}\geq\left[  \frac{\left(  \frac{d_{2}}{\lambda_{k}d_{1}}\right)
^{b_{k}}\left(  1-\frac{d_{2}}{\lambda_{k}d_{1}}\right)  ^{1-b_{k}}}
{b_{k}^{b_{k}}\left(  1-b_{k}\right)  ^{1-b_{k}}}\right]  ^{\frac{1}{b_{k}
}\cdot\frac{s+1}{qr-\left(  p-1\right)  \left(  s+1\right)  }}
\]
where
\[
b_{k}=\frac{s+1-\lambda_{k}r}{q-\lambda_{k}\left(  p-1\right)  }\in\left(
0,1\right)  .
\]
Since
\begin{align*}
&  \lim_{k\rightarrow\infty}\left[  \frac{\left(  \frac{d_{2}}{\lambda
_{k}d_{1}}\right)  ^{b_{k}}\left(  1-\frac{d_{2}}{\lambda_{k}d_{1}}\right)
^{1-b_{k}}}{b_{k}^{b_{k}}\left(  1-b_{k}\right)  ^{1-b_{k}}}\right]
^{\frac{1}{b_{k}}}=\lim_{k\rightarrow\infty}\left(  \frac{1-\frac{d_{2}
}{\lambda_{k}d_{1}}}{1-b_{k}}\right)  ^{1-b_{k}}\\
=  &  \lim_{k\rightarrow\infty}\left(  \frac{qr-\left(  p-1\right)  \left(
s+1\right)  }{\lambda_{k}\left(  q-\lambda_{k}\left(  p-1\right)  \right)
\left(  q-\lambda\left(  p-1\right)  \right)  }\right)  ^{1-b_{k}}=1,
\end{align*}
we again have $\underline{u}\geq1$. Optimal bound for $v$ follows from
$\left(  \ref{Equation simple inequality 1}\right)  $.
\end{proof}

\begin{remark}
The admissible set $\Lambda_{2}$ is a nonempty interval with right end point
$\frac{s+1}{r}$, from the property of function $f_{0}$, the set $f_{0}\left(
\Lambda_{2}\right)  $ is a nonempty interval of the form $\left(
0,k_{2}\right]  $ for some constant $k_{2}$ depending on $p,q,r,s$. If
$\frac{p-1}{r}\geq1$, it is easy to check
\[
\Lambda_{2}=\left[  1,\frac{s+1}{r}\right)  ,
\]
and
\[
f_{0}\left(  \Lambda_{2}\right)  \supset\left(  0,\frac{s+1-r}{q-\left(
p-1\right)  }\right]  .
\]
And if $\frac{p-1}{r}<1$,
\[
\Lambda_{2}=\left[  \max\left\{  1,\frac{\left(  s+1\right)  -q}{r-\left(
p-1\right)  }\right\}  ,\frac{s+1}{r}\right)
\]
so if in addition, $\frac{\left(  s+1\right)  -q}{r-\left(  p-1\right)  }<1$,
we have
\[
\Lambda_{2}=\left[  1,\frac{s+1}{r}\right)  ,\quad f_{0}\left(  \Lambda
_{2}\right)  \supset\left(  0,\frac{s+1-r}{q-\left(  p-1\right)  }\right]  ,
\]
and if in addition, $\frac{\left(  s+1\right)  -q}{r-\left(  p-1\right)  }
\geq1$, we have
\[
\Lambda_{2}=\left[  \frac{\left(  s+1\right)  -q}{r-\left(  p-1\right)
},\frac{s+1}{r}\right)  ,\quad f_{0}\left(  \Lambda_{2}\right)  \supset\left(
0,\frac{s+1-q}{r-\left(  p-1\right)  }\right]  .
\]
When $\left(  p,q,r,s\right)  =\left(  2,4,2,4\right)  $, $\Lambda_{2}=\left(
1,\frac{5}{2}\right)  $ and $f_{0}\left(  \Lambda_{2}\right)  =\left(
0,11-4\sqrt{6}\right)  $.
\end{remark}

Combining the optimal bounds in Theorems \ref{Theorem upper bounds} and
\ref{Theorem lower bounds}, we have

\begin{theorem}
Let
\[
s+1>\max\left\{  q,r\right\}  .
\]
Then $u\equiv1,v\equiv1$ is the only solution whenever
\[
\frac{d_{2}}{d_{1}}\in f_{0}\left(  \Lambda_{1}\right)  \cap f_{0}\left(
\Lambda_{2}\right)  .
\]

\end{theorem}

\begin{remark}
In general, $f_{0}\left(  \Lambda_{1}\right)  =\left(  0,k_{1}\right]  $ and
$f_{0}\left(  \Lambda_{2}\right)  =\left(  0,k_{2}\right]  $, hence
\[
f_{0}\left(  \Lambda_{1}\right)  \cap f_{0}\left(  \Lambda_{2}\right)
=\left(  0,k\right]
\]
where $k=\min\left\{  k_{1},k_{2}\right\}  $. When $\left(  p,q,r,s\right)
=\left(  2,4,2,4\right)  $, we have $k=1$, hence $u\equiv1,v\equiv1$ is the
only solution when $d_{2}\leq d_{1}$.
\end{remark}

Now we extend our optimal bounds to the case $\sigma>0$.

\begin{theorem}
\label{Theorem upper bounds positive}Assume $\sigma>0$ and $q<s+1$. If
\[
\frac{d_{2}}{d_{1}}\in f_{\sigma}\left(  \Lambda_{3}\right)  ,
\]
where
\[
\Lambda_{3}=\left(  0,1\right]  \cap\left(  0,\frac{s+1-q}{r-\left(
p-1\right)  }\right]  \cap\left(  0,\frac{s+1}{r+1}\right]  ,
\]
then $u\leq u^{\ast},v\leq v^{\ast}$.
\end{theorem}

\begin{proof}
Let $\lambda\in\tilde{\Lambda}_{3}\cap\left(  0,\frac{d_{2}}{d_{1}}\right)  $
where
\[
\tilde{\Lambda}_{3}=\left(  0,1\right]  \cap\left(  0,\frac{s+1-q}{r-\left(
p-1\right)  }\right)  \cap\left(  0,\frac{s+1}{r+1}\right)  .
\]
Since $0<\lambda\leq1$, we have from Lemma \ref{Lemma local maximum}, at any
point $x^{\ast}\in\overline{\Omega}$ where $\frac{u}{v^{\lambda}}$ achieves
its maximum,
\[
1-\frac{u^{p-1}}{v^{q}}-\frac{\lambda d_{1}}{d_{2}}+\frac{\lambda d_{1}}
{d_{2}}\frac{u^{r}}{v^{s+1}}-\sigma u^{-1}\leq0.
\]
Since $\lambda\in\left(  0,\frac{s+1-q}{r-\left(  p-1\right)  }\right)
\cap\left(  0,\frac{s+1}{r+1}\right)  $, we also have
\[
a=\frac{q-\lambda\left(  p-1\right)  }{s+1-\lambda r}\in\left(  0,1\right)
\text{ and }a_{0}=\frac{\lambda}{s+1-\lambda r}\in\left(  0,1\right)  .
\]
Let $\delta$ be any given number in $\left(  0,1\right)  $, we have from
Young's inequality
\begin{align*}
&  \frac{u^{p-1}}{v^{q}}=\left(  \frac{\delta}{a}\frac{\lambda d_{1}}{d_{2}
}\frac{u^{r}}{v^{s+1}}\right)  ^{a}\left(  \left(  \frac{\delta}{a}
\frac{\lambda d_{1}}{d_{2}}\right)  ^{-\frac{a}{1-a}}\left(  \frac
{u}{v^{\lambda}}\right)  ^{-\frac{qr-\left(  p-1\right)  \left(  s+1\right)
}{s+1-\lambda r-\left(  q-\lambda\left(  p-1\right)  \right)  }}\right)
^{1-a}\\
\leq &  \delta\frac{\lambda d_{1}}{d_{2}}\frac{u^{r}}{v^{s+1}}+\frac
{1-a}{\left(  \frac{\delta}{a}\frac{\lambda d_{1}}{d_{2}}\right)  ^{\frac
{a}{1-a}}}\left(  \frac{u}{v^{\lambda}}\right)  ^{-\frac{qr-\left(
p-1\right)  \left(  s+1\right)  }{s+1-\lambda r-\left(  q-\lambda\left(
p-1\right)  \right)  }},
\end{align*}
and
\begin{align*}
&  \sigma u^{-1}=\left(  \frac{1-\delta}{a_{0}}\frac{\lambda d_{1}}{d_{2}
}\frac{u^{r}}{v^{s+1}}\right)  ^{a_{0}}\left(  \sigma^{\frac{1}{1-a_{0}}
}\left(  \frac{1-\delta}{a_{0}}\frac{\lambda d_{1}}{d_{2}}\right)
^{-\frac{a_{0}}{1-a_{0}}}\left(  \frac{u}{v^{\lambda}}\right)  ^{-\frac
{s+1}{s+1-\lambda r-\lambda}}\right)  ^{1-a_{0}}\\
\leq &  \left(  1-\delta\right)  \frac{\lambda d_{1}}{d_{2}}\frac{u^{r}
}{v^{s+1}}+\frac{\sigma^{\frac{1}{1-a_{0}}}\left(  1-a_{0}\right)  }{\left(
\frac{1-\delta}{a_{0}}\frac{\lambda d_{1}}{d_{2}}\right)  ^{\frac{a_{0}
}{1-a_{0}}}}\left(  \frac{u}{v^{\lambda}}\right)  ^{-\frac{s+1}{s+1-\lambda
r-\lambda}}.
\end{align*}
Combining the three inequalities above and applying Lemma
\ref{Lemma from u/v^lambda to u,v}, we deduce
\begin{equation}
1-\frac{\lambda d_{1}}{d_{2}}\leq\frac{1-a}{\left(  \frac{\delta}{a}
\frac{\lambda d_{1}}{d_{2}}\right)  ^{\frac{a}{1-a}}}\bar{u}^{-\frac
{qr-\left(  p-1\right)  \left(  s+1\right)  }{s+1}\frac{1}{1-a}}+\frac
{\sigma^{\frac{1}{1-a_{0}}}\left(  1-a_{0}\right)  }{\left(  \frac{1-\delta
}{a_{0}}\frac{\lambda d_{1}}{d_{2}}\right)  ^{\frac{a_{0}}{1-a_{0}}}}\bar
{u}^{-\frac{1}{1-a_{0}}}\label{Equation inequality for u positive}%
\end{equation}
which yields an upper bound for $\bar{u}$ since
\[
1-\frac{\lambda d_{1}}{d_{2}}>0.
\]
Upper bound for $\bar{v}$ follows from $\left(
\ref{Equation simple inequality 1}\right)  $.

Next, if $\frac{d_{2}}{d_{1}}\in f_{\sigma}\left(  \Lambda_{3}\right)  $, then
there exists $\lambda\in\Lambda_{3}$ such that $\frac{d_{2}}{d_{1}}=f_{\sigma
}\left(  \lambda\right)  $. We first assume $\lambda\in\tilde{\Lambda}_{3}$,
then $\left(  \ref{Equation inequality for u positive}\right)  $ holds for any
$\delta\in\left(  0,1\right)  $. Let
\[
\delta=a\frac{d_{2}}{\lambda d_{1}}\left(  u^{\ast}\right)  ^{-\frac
{qr-\left(  p-1\right)  \left(  s+1\right)  }{s+1}}\equiv1-a_{0}\sigma
\frac{d_{2}}{\lambda d_{1}}\left(  u^{\ast}\right)  ^{-1}\in\left(
0,1\right)  ,
\]
then $\left(  \ref{Equation inequality for u positive}\right)  $ becomes
\begin{equation}
1-\frac{\lambda d_{1}}{d_{2}}\leq\left(  1-a\right)  \left(  \frac{\bar{u}
}{\left(  u^{\ast}\right)  ^{a}}\right)  ^{-\frac{qr-\left(  p-1\right)
\left(  s+1\right)  }{s+1}\frac{1}{1-a}}+\sigma\left(  1-a_{0}\right)  \left(
\frac{\bar{u}}{\left(  u^{\ast}\right)  ^{a_{0}}}\right)  ^{-\frac{1}{1-a_{0}
}}.\label{Equation inequality for u positive
upper}%
\end{equation}
We observe that if $\bar{u}=u^{\ast}$, then
\begin{align*}
&  \left(  1-a\right)  \left(  \frac{u^{\ast}}{\left(  u^{\ast}\right)  ^{a}
}\right)  ^{-\frac{qr-\left(  p-1\right)  \left(  s+1\right)  }{s+1}\frac
{1}{1-a}}+\sigma\left(  1-a_{0}\right)  \left(  \frac{u^{\ast}}{\left(
u^{\ast}\right)  ^{a_{0}}}\right)  ^{-\frac{1}{1-a_{0}}}\\
=  &  \left(  1-a\right)  \left(  u^{\ast}\right)  ^{-\frac{qr-\left(
p-1\right)  \left(  s+1\right)  }{s+1}}+\left(  1-a_{0}\right)  \sigma\left(
u^{\ast}\right)  ^{-1}\\
=  &  \left(  u^{\ast}\right)  ^{-\frac{qr-\left(  p-1\right)  \left(
s+1\right)  }{s+1}}+\sigma\left(  u^{\ast}\right)  ^{-1}-\left(  a\left(
u^{\ast}\right)  ^{-\frac{qr-\left(  p-1\right)  \left(  s+1\right)  }{s+1}
}+a_{0}\sigma\left(  u^{\ast}\right)  ^{-1}\right) \\
=  &  1-\frac{\lambda d_{1}}{d_{2}}.
\end{align*}
Since the right hand side of $\left(
\ref{Equation inequality for u positive upper}\right)  $ is monotone
decreasing in $\bar{u}$, the equality above implies $\bar{u}\leq u^{\ast}$.
And from $\left(  \ref{Equation simple inequality 1}\right)  $, we have
$\bar{v}\leq\left(  u^{\ast}\right)  ^{\frac{r}{s+1}}=v^{\ast}$. If
$\lambda\not \in \tilde{\Lambda}_{3}$, then $\lambda=\frac{s+1-q}{r-\left(
p-1\right)  }$ or $\frac{s+1}{r+1}$, and there exist $\lambda_{k}\in
\tilde{\Lambda}_{3}$ such that $\lambda_{k}<\lambda$ and
\[
\lim_{k\rightarrow\infty}\lambda_{k}=\lambda.
\]
Let
\[
a_{k}=\frac{q-\lambda_{k}\left(  p-1\right)  }{s+1-\lambda_{k}r}\in\left(
0,1\right)  \text{ and }a_{0k}=\frac{\lambda_{k}}{s+1-\lambda_{k}r}\in\left(
0,1\right)  ,
\]
we have
\[
1-\frac{\lambda_{k}d_{1}}{d_{2}}\leq\left(  1-a_{k}\right)  \left(  \frac
{\bar{u}}{\left(  u^{\ast}\right)  ^{a_{k}}}\right)  ^{-\frac{qr-\left(
p-1\right)  \left(  s+1\right)  }{s+1}\frac{1}{1-a_{k}}}+\sigma\left(
1-a_{0k}\right)  \left(  \frac{\bar{u}}{\left(  u^{\ast}\right)  ^{a_{0k}}
}\right)  ^{-\frac{1}{1-a_{0k}}}.
\]
For each $k$, let $u_{k}$ be the unique positive number such that
\begin{align}
&  1-\frac{\lambda_{k}d_{1}}{d_{2}}\label{Equation uk}\\
=  &  \left(  1-a_{k}\right)  \left(  \frac{u_{k}}{\left(  u^{\ast}\right)
^{a_{k}}}\right)  ^{-\frac{qr-\left(  p-1\right)  \left(  s+1\right)  }
{s+1}\frac{1}{1-a_{k}}}+\sigma\left(  1-a_{0k}\right)  \left(  \frac{u_{k}
}{\left(  u^{\ast}\right)  ^{a_{0k}}}\right)  ^{-\frac{1}{1-a_{0k}}},\nonumber
\end{align}
then we have $\bar{u}\leq u_{k}$. We claim
\[
\lim_{k\rightarrow\infty}u_{k}=u^{\ast}.
\]
First, since $\left(  u^{\ast},v^{\ast}\right)  $ is a solution, $\bar{u}\leq
u_{k}$ implies that for each $k$, $u_{k}\geq u^{\ast}$. If
\[
\lambda=\frac{s+1-q}{r-\left(  p-1\right)  }<\frac{s+1}{r+1},
\]
then direct calculation yields
\[
1-\frac{\lambda_{k}d_{1}}{d_{2}}>1-\frac{\lambda d_{1}}{d_{2}}=\sigma\left(
u^{\ast}\right)  ^{-1}\frac{s+1-\lambda\left(  r+1\right)  }{s+1-\lambda
r}>0.
\]
It is easy to check
\[
\lim_{k\rightarrow\infty}a_{k}=1,\lim_{k\rightarrow\infty}a_{0k}=a_{0}
\in\left(  0,1\right)  .
\]
Hence
\[
\left(  1-a_{k}\right)  \left(  \frac{u_{k}}{\left(  u^{\ast}\right)  ^{a_{k}
}}\right)  ^{-\frac{qr-\left(  p-1\right)  \left(  s+1\right)  }{s+1}\frac
{1}{1-a_{k}}}\leq\left(  1-a_{k}\right)  \left(  u^{\ast}\right)
^{-\frac{qr-\left(  p-1\right)  \left(  s+1\right)  }{s+1}}
\]
which converges to zero as $k\rightarrow\infty$. So we have $\lim
_{k\rightarrow\infty}u_{k}$ exists, and
\[
1-\frac{\lambda d_{1}}{d_{2}}=\sigma\left(  1-a_{0}\right)  \left(  \frac
{\lim_{k\rightarrow\infty}u_{k}}{\left(  u^{\ast}\right)  ^{a_{0}}}\right)
^{-\frac{1}{1-a_{0}}},
\]
hence $\lim_{k\rightarrow\infty}u_{k}=u^{\ast}$. If
\[
\lambda=\frac{s+1}{r+1}<\frac{s+1-q}{r-\left(  p-1\right)  },
\]
we have
\[
\lim_{k\rightarrow\infty}a_{k}=a\in\left(  0,1\right)  ,\quad\lim
_{k\rightarrow\infty}a_{0k}=1.
\]
We can obtain $\lim_{k\rightarrow\infty}u_{k}=u^{\ast}$ similarly. Finally, we
consider the case
\[
\lambda=\frac{s+1}{r+1}=\frac{s+1-q}{r-\left(  p-1\right)  },
\]
then
\[
\lim_{k\rightarrow\infty}a_{k}=1,\quad\lim_{k\rightarrow\infty}a_{0k}=1
\]
and
\[
\frac{\lambda d_{1}}{d_{2}}=1,
\]
dividing equation $\left(  \ref{Equation uk}\right)  $ by $\lambda-\lambda
_{k}$, we have
\begin{align*}
&  \frac{d_{1}}{d_{2}}\\
=  &  \frac{1-a_{k}}{\left(  \lambda-\lambda_{k}\right)  \left(  \frac
{\delta_{k}}{a_{k}}\frac{\lambda_{k}d_{1}}{d_{2}}\right)  ^{\frac{a_{k}
}{1-a_{k}}}}u_{k}^{-\frac{qr-\left(  p-1\right)  \left(  s+1\right)  }
{s+1}\frac{1}{1-a_{k}}}+\frac{\sigma^{\frac{1}{1-a_{0k}}}\left(
1-a_{0k}\right)  }{\left(  \lambda-\lambda_{k}\right)  \left(  \frac
{1-\delta_{k}}{a_{0k}}\frac{\lambda_{k}d_{1}}{d_{2}}\right)  ^{\frac{a_{0k}
}{1-a_{0k}}}}u_{k}^{-\frac{1}{1-a_{0k}}}\\
=  &  \frac{qr-\left(  p-1\right)  \left(  s+1\right)  }{\left(  s+1-\lambda
r\right)  \left(  s+1-\lambda_{k}r\right)  }\left(  \frac{u_{k}}{\left(
u^{\ast}\right)  ^{a_{k}}}\right)  ^{-\frac{qr-\left(  p-1\right)  \left(
s+1\right)  }{s+1}\frac{1}{1-a_{k}}}\\
&  +\frac{\sigma\left(  s+1\right)  }{\left(  s+1-\lambda r\right)  \left(
s+1-\lambda_{k}r\right)  }\left(  \frac{u_{k}}{\left(  u^{\ast}\right)
^{a_{0k}}}\right)  ^{-\frac{1}{1-a_{0k}}}.
\end{align*}
If
\[
\liminf_{k\rightarrow\infty}u_{k}>u^{\ast},
\]
then both terms in the right hand side tend to zero, which is a contradiction;
If
\[
\limsup_{k\rightarrow\infty}u_{k}<u^{\ast},
\]
then both terms in the right hand side tend to infinity, which is a also
contradiction. Hence, we have $\lim_{k\rightarrow\infty}u_{k}=u^{\ast}$.
\end{proof}

Next, we extend optimal lower bound to the case $\sigma>0$.

\begin{theorem}
\label{Theorem lower bounds positive}If $s+1>r$ and $\frac{d_{2}}{d_{1}}\in
f_{\sigma}\left(  \Lambda_{4}\right)  $, where
\[
\Lambda_{4}=\left\{  \lambda\in\left[  1,\infty\right)  :\frac{s+1}{r+1}
\leq\lambda<\frac{s+1}{r},\frac{s+1-\lambda r}{q-\lambda\left(  p-1\right)
}\leq1\right\}
\]
then we have $u\geq u^{\ast}$, $v\geq v^{\ast}$.
\end{theorem}

\begin{proof}
Let $\lambda\in\tilde{\Lambda}_{4}\cap\left(  \frac{d_{2}}{d_{1}}
,\infty\right)  $ where
\[
\tilde{\Lambda}_{4}=\left\{  \lambda\in\left[  1,\infty\right)  :\frac
{s+1}{r+1}<\lambda<\frac{s+1}{r},\frac{s+1-\lambda r}{q-\lambda\left(
p-1\right)  }<1\right\}  .
\]
Since $\lambda\geq1$, we have from Lemma \ref{Lemma local minimum}, at any
point $x^{\ast}\in\overline{\Omega}$ where $\frac{u}{v^{\lambda}}$ achieves
its minimum,
\[
1-\frac{u^{p-1}}{v^{q}}-\sigma u^{-1}-\frac{\lambda d_{1}}{d_{2}}\left(
1-\frac{u^{r}}{v^{s+1}}\right)  \geq0.
\]
Under our assumptions, we also have
\[
b=\frac{s+1-\lambda r}{q-\lambda\left(  p-1\right)  }\in\left(  0,1\right)
,\quad b_{0}=\frac{s+1-\lambda r}{\lambda}\in\left(  0,1\right)  .
\]
Applying Young's inequality, we have
\begin{align*}
&  \frac{\lambda d_{1}}{d_{2}}-1\leq\frac{\lambda d_{1}}{d_{2}}\frac{u^{r}
}{v^{s+1}}-\frac{u^{p-1}}{v^{q}}-\sigma u^{-1}\\
=  &  \left[  \frac{1}{b}\frac{u^{p-1}}{v^{q}}\right]  ^{b}\left(  \left(
\delta\frac{\lambda d_{1}}{d_{2}}\right)  ^{\frac{1}{1-b}}b^{\frac{b}{1-b}
}\left(  \frac{u}{v^{\lambda}}\right)  ^{\frac{qr-\left(  p-1\right)  \left(
s+1\right)  }{q-\lambda\left(  p-1\right)  -\left(  s+1-\lambda r\right)  }
}\right)  ^{1-b}\\
&  +\left[  \frac{1}{b_{0}}\sigma u^{-1}\right]  ^{b_{0}}\left(
\sigma^{-\frac{b_{0}}{1-b_{0}}}\left(  \left(  1-\delta\right)  \frac{\lambda
d_{1}}{d_{2}}\right)  ^{\frac{1}{1-b_{0}}}b_{0}^{\frac{b_{0}}{1-b_{0}}}\left(
\frac{u}{v^{\lambda}}\right)  ^{\frac{s+1}{\lambda-\left(  s+1-\lambda
r\right)  }}\right)  ^{1-b_{0}}\\
&  -\frac{u^{p-1}}{v^{q}}-\sigma u^{-1}\\
\leq &  \left(  1-b\right)  \left(  \delta\frac{\lambda d_{1}}{d_{2}}\right)
^{\frac{1}{1-b}}b^{\frac{b}{1-b}}\left(  \frac{u}{v^{\lambda}}\right)
^{\frac{qr-\left(  p-1\right)  \left(  s+1\right)  }{q-\lambda\left(
p-1\right)  -\left(  s+1-\lambda r\right)  }}\\
&  +\left(  1-b_{0}\right)  \sigma^{-\frac{b_{0}}{1-b_{0}}}\left(  \left(
1-\delta\right)  \frac{\lambda d_{1}}{d_{2}}\right)  ^{\frac{1}{1-b_{0}}}
b_{0}^{\frac{b_{0}}{1-b_{0}}}\left(  \frac{u}{v^{\lambda}}\right)
^{\frac{s+1}{\lambda-\left(  s+1-\lambda r\right)  }}\\
\leq &  \left(  1-b\right)  \left(  \delta\frac{\lambda d_{1}}{d_{2}}\right)
^{\frac{1}{1-b}}b^{\frac{b}{1-b}}\underline{u}^{\frac{qr-\left(  p-1\right)
\left(  s+1\right)  }{s+1}\frac{b}{1-b}}\\
&  +\left(  1-b_{0}\right)  \sigma^{-\frac{b_{0}}{1-b_{0}}}\left(  \left(
1-\delta\right)  \frac{\lambda d_{1}}{d_{2}}\right)  ^{\frac{1}{1-b_{0}}}
b_{0}^{\frac{b_{0}}{1-b_{0}}}\underline{u}^{\frac{b_{0}}{1-b_{0}}}%
\end{align*}
which yields\ a lower bound of $\underline{u}$ since
\[
\frac{\lambda d_{1}}{d_{2}}-1>0.
\]
Lower bound of $\underline{v}$ follows from $\left(
\ref{Equation simple inequality 1}\right)  $. When
\[
\frac{d_{2}}{d_{1}}=f_{\sigma}\left(  \lambda\right)
\]
for some $\lambda\in\Lambda_{4}$, $u\geq u^{\ast}$, $v\geq v^{\ast}$ can be
proved similarly as in the proof of Theorem
\ref{Theorem upper bounds positive} with the same choice of $\delta$.
\end{proof}

Combining Theorems \ref{Theorem upper bounds positive} and
\ref{Theorem lower bounds positive}, we have

\begin{corollary}
\label{Corollary optimal when sigma>0}Assume $\sigma>0$ and $\max\left\{
q,r\right\}  <s+1$. Then we have $u\equiv u^{\ast},v\equiv v^{\ast}$ if
\[
\frac{d_{2}}{d_{1}}\in f_{\sigma}\left(  \Lambda_{3}\right)  \cap f_{\sigma
}\left(  \Lambda_{4}\right)  .
\]

\end{corollary}

\begin{remark}
It is easy to see that $f_{\sigma}\left(  \Lambda_{3}\right)  =\left(
0,k_{1}\right]  $ and $f_{\sigma}\left(  \Lambda_{4}\right)  =\left(
0,k_{2}\right]  $ for some $k_{1},k_{2}>0$. Hence, $f_{\sigma}\left(
\Lambda_{3}\right)  \cap f_{\sigma}\left(  \Lambda_{4}\right)  =\left(
0,k\right]  $ with $k=\min\left\{  k_{1},k_{2}\right\}  >0$.
\end{remark}

When $\sigma>0$, there are simple lower bounds
\[
u>\sigma,\quad v>\sigma^{\frac{r}{s+1}}.
\]
In Theorem \ref{Theorem upper bounds positive}, we haven't made use of such
bounds. Actually, these lower bounds will help to relax the conditions in
Theorem \ref{Theorem upper bounds positive}. However, we no longer obtain
optimal bounds.

\begin{theorem}
\label{Theorem upper bounds positive copy(1)}Assume $\sigma>0$ and $p-1<r$.
Then
\[
u\leq c\left(  1+\left(  \frac{d_{2}}{d_{1}}\right)  ^{\gamma}\right)  ,\quad
v\leq c\left(  1+\left(  \frac{d_{2}}{d_{1}}\right)  ^{\frac{s+1}{r}\gamma
}\right)  ,
\]
where $c,\gamma$ are positive constants independent of $d_{1},d_{2}$.
\end{theorem}

\begin{proof}
Let $\varepsilon\geq0$ be such that
\[
0<\frac{p-1}{r}<\frac{q-\varepsilon}{s+1}<1.
\]
Let $\lambda\in\left(  0,\frac{d_{2}}{d_{1}}\right)  \cap\Lambda_{5}$ where
\[
\Lambda_{5}=\left\{  \lambda\in\left(  0,1\right]  :\lambda<\min\left\{
\frac{s+1-\left(  q-\varepsilon\right)  }{r-\left(  p-1\right)  },\frac
{s+1}{r+1}\right\}  \right\}  .
\]
Since $0<\lambda\leq1$, we have from Lemma \ref{Lemma local maximum}, at any
point $x^{\ast}\in\overline{\Omega}$ where $\frac{u}{v^{\lambda}}$ achieves
its maximum,
\[
1-\frac{u^{p-1}}{v^{q}}-\frac{\lambda d_{1}}{d_{2}}+\frac{\lambda d_{1}}
{d_{2}}\frac{u^{r}}{v^{s+1}}-\sigma u^{-1}\leq0.
\]
Since
\[
a_{\varepsilon}=\frac{q-\varepsilon-\lambda\left(  p-1\right)  }{s+1-\lambda
r}\in\left(  0,1\right)  \text{ and }a_{0}=\frac{\lambda}{s+1-\lambda r}
\in\left(  0,1\right)  ,
\]
for any $\delta\in\left(  0,1\right)  $, we have from Young's inequality,
\begin{align*}
&  \frac{u^{p-1}}{v^{q}}<\sigma^{\frac{r\varepsilon}{s+1}}\frac{u^{p-1}
}{v^{q-\varepsilon}}\\
=  &  \left(  \frac{\delta}{a_{\varepsilon}}\frac{\lambda d_{1}}{d_{2}}
\frac{u^{r}}{v^{s+1}}\right)  ^{a_{\varepsilon}}\left(  \sigma^{\frac{r}
{s+1}\frac{\varepsilon}{1-a_{\varepsilon}}}\left(  \frac{\delta}
{a_{\varepsilon}}\frac{\lambda d_{1}}{d_{2}}\right)  ^{-\frac{a_{\varepsilon}
}{1-a_{\varepsilon}}}\left(  \frac{u}{v^{\lambda}}\right)  ^{-\frac{\left(
q-\varepsilon\right)  r-\left(  p-1\right)  \left(  s+1\right)  }{s+1-\lambda
r-\left(  q-\varepsilon-\lambda\left(  p-1\right)  \right)  }}\right)
^{1-a_{\varepsilon}}\\
\leq &  \delta\frac{\lambda d_{1}}{d_{2}}\frac{u^{r}}{v^{s+1}}+\frac{\left(
1-a_{\varepsilon}\right)  \sigma^{\frac{r}{s+1}\frac{\varepsilon
}{1-a_{\varepsilon}}}}{\left(  \frac{\delta}{a_{\varepsilon}}\frac{\lambda
d_{1}}{d_{2}}\right)  ^{\frac{a_{\varepsilon}}{1-a_{\varepsilon}}}}\left(
\frac{u}{v^{\lambda}}\right)  ^{-\frac{\left(  q-\varepsilon\right)  r-\left(
p-1\right)  \left(  s+1\right)  }{s+1-\lambda r-\left(  q-\varepsilon
-\lambda\left(  p-1\right)  \right)  }}%
\end{align*}
and
\begin{align*}
&  \sigma u^{-1}=\left(  \frac{1-\delta}{a_{0}}\frac{\lambda d_{1}}{d_{2}
}\frac{u^{r}}{v^{s+1}}\right)  ^{a_{0}}\left(  \sigma^{\frac{1}{1-a_{0}}
}\left(  \frac{1-\delta}{a_{0}}\frac{\lambda d_{1}}{d_{2}}\right)
^{-\frac{a_{0}}{1-a_{0}}}\left(  \frac{u}{v^{\lambda}}\right)  ^{-\frac
{s+1}{s+1-\lambda r-\lambda}}\right)  ^{1-a_{0}}\\
\leq &  \left(  1-\delta\right)  \frac{\lambda d_{1}}{d_{2}}\frac{u^{r}
}{v^{s+1}}+\frac{\sigma^{\frac{1}{1-a_{0}}}\left(  1-a_{0}\right)  }{\left(
\frac{1-\delta}{a_{0}}\frac{\lambda d_{1}}{d_{2}}\right)  ^{\frac{a_{0}
}{1-a_{0}}}}\left(  \frac{u}{v^{\lambda}}\right)  ^{-\frac{s+1}{s+1-\lambda
r-\lambda}}.
\end{align*}
Combining the above inequalities and applying Lemma
\ref{Lemma from u/v^lambda to u,v}, we have
\[
1-\frac{\lambda d_{1}}{d_{2}}\leq\frac{\left(  1-a_{\varepsilon}\right)
\sigma^{\frac{r}{s+1}\frac{\varepsilon}{1-a_{\varepsilon}}}}{\left(
\frac{\delta}{a_{\varepsilon}}\frac{\lambda d_{1}}{d_{2}}\right)
^{\frac{a_{\varepsilon}}{1-a_{\varepsilon}}}}\left(  \bar{u}\right)
^{-\frac{\left(  q-\varepsilon\right)  r-\left(  p-1\right)  \left(
s+1\right)  }{s+1}\cdot\frac{1}{1-a_{\varepsilon}}}+\frac{\sigma^{\frac
{1}{1-a_{0}}}\left(  1-a_{0}\right)  }{\left(  \frac{1-\delta}{a_{0}}
\frac{\lambda d_{1}}{d_{2}}\right)  ^{\frac{a_{0}}{1-a_{0}}}}\bar{u}
^{-\frac{1}{1-a_{0}}},
\]
which yields an upper bound of $\bar{u}$ since
\[
1-\frac{\lambda d_{1}}{d_{2}}>0.
\]
Upper bound for $\bar{v}$ follows from $\left(
\ref{Equation simple inequality 1}\right)  $.
\end{proof}

Theorem \ref{Theorem main} and Theorem \ref{Theorem main 1} are proved by
combining all the results we have so far in this section.

In the remaining part of this section, we consider $\left(
\ref{Equation stationary}\right)  $ with common sources, i.e., we require
$p=r$, $q=s$. The assumption
\[
\frac{p-1}{r}<\frac{q}{s+1}
\]
is then reduced to $s+1>r$. And $\Lambda_{i}$, $i=1,2,3,4$ have simpler forms
\[%
\begin{array}
[c]{ll}%
\Lambda_{1}=\left(  0,1\right]  , & \Lambda_{2}=\left[  1,\frac{s+1}
{r}\right)  ,\\
\Lambda_{3}=\left(  0,1\right]  \cap\left(  0,\frac{s+1}{r+1}\right]  , &
\Lambda_{4}=\left[  1,\frac{s+1}{r}\right)  \cap\left[  \frac{s+1}{r+1}
,\frac{s+1}{r}\right)  .
\end{array}
\]
For any $\lambda\in\left[  0,\frac{s+1}{r}\right]  $,
\[
f_{\sigma}\left(  \lambda\right)  =\frac{\lambda\left(  s+1-\lambda r\right)
}{\left(  s-\lambda r\right)  \left(  u^{\ast}\right)  ^{\frac{r}{s+1}
-1}+\lambda}.
\]
Obviously,
\[
f_{\sigma}\left(  \frac{s}{r}\right)  \equiv f_{0}\left(  \frac{s}{r}\right)
=1.
\]
Since $u^{\ast}$ is monotone increasing in $\sigma$, we have for any fixed
$\lambda\in\left(  0,\frac{s}{r}\right)  $, $f_{\sigma}\left(  \lambda\right)
$ is monotone increasing in $\sigma$; for any $\lambda\in\left(  \frac{s}
{r},\frac{s+1}{r}\right)  $, $f_{\sigma}\left(  \lambda\right)  $ is monotone
decreasing in $\sigma$. Furthermore, since $\lim_{\sigma\rightarrow\infty
}u^{\ast}=\infty$, for any $\lambda\in\left(  0,\frac{s+1}{r}\right)  $,
\[
\lim_{\sigma\rightarrow\infty}f_{\sigma}\left(  \lambda\right)  =s+1-\lambda
r.
\]

Let $k_{\sigma}$ be such that for $\sigma=0$, $\left(  0,k_{0}\right]
=f_{0}\left(  \Lambda_{1}\right)  \cap f_{0}\left(  \Lambda_{2}\right)  $ and
for $\sigma>0$, $\left(  0,k_{\sigma}\right]  =f_{\sigma}\left(  \Lambda
_{3}\right)  \cap f_{\sigma}\left(  \Lambda_{4}\right)  $. Then for any
$\sigma\geq0$, the system has no nontrivial solution whenever
\[
\frac{d_{2}}{d_{1}}\leq k_{\sigma}.
\]

The following theorem describes the dependence of $k_{\sigma}$ on $\sigma$.

\begin{theorem}
\label{Theorem example}Let $\sigma\geq0$.

\begin{enumerate}
\item[(i)] $k_{0}=1$.

\item[(ii)] If $s=r$, then for any $\sigma>0$, $k_{\sigma}=1$.

\item[(iii)] If $s>r$, then $k_{\sigma}$ is monotone increasing for $\sigma
\in\left(  0,\infty\right)  $, and
\[
\lim_{\sigma\rightarrow0}k_{\sigma}=k_{0}=1,\quad\lim_{\sigma\rightarrow
\infty}k_{\sigma}=\frac{s+1}{r+1}.
\]

\item[(iv)] If $s<r$, then $k_{\sigma}$ is monotone decreasing for $\sigma
\in\left(  0,\infty\right)  $, and
\[
\lim_{\sigma\rightarrow0}k_{\sigma}=k_{0}=1,\quad\lim_{\sigma\rightarrow
\infty}k_{\sigma}=s+1-r.
\]

\end{enumerate}
\end{theorem}

\begin{proof}
(i). We have
\[
\Lambda_{1}=\left(  0,1\right]  ,\Lambda_{2}=\left[  1,\frac{s+1}{r}\right)
.
\]
Since $f_{0}$ has only one critical point in $\left(  0,\frac{s+1}{r}\right)
$, we have
\[
\left(  0,k_{0}\right]  =f_{0}\left(  \Lambda_{1}\right)  \cap f_{0}\left(
\Lambda_{2}\right)  =\left(  0,f_{0}\left(  1\right)  \right]  =\left(
0,1\right]  ,
\]
hence $k_{0}=1$.

(ii). If $s=r$, we have
\[
\Lambda_{3}=\left(  0,1\right]  \cap\left(  0,\frac{s+1}{r+1}\right]  =\left(
0,1\right]  ,
\]
and
\[
\Lambda_{4}=\left[  1,\frac{s+1}{r}\right)  \cap\left[  \frac{s+1}{r+1}
,\frac{s+1}{r}\right)  =\left[  1,\frac{s+1}{r}\right)  .
\]
From the property of $f_{\sigma}$, we have $k_{\sigma}=f_{\sigma}\left(
1\right)  =1$.

(iii). $s>r$ implies
\[
0<1<\frac{s+1}{r+1}<\frac{s}{r}<\frac{s+1}{r}.
\]
Hence
\[
\Lambda_{3}=\left(  0,1\right]  \cap\left(  0,\frac{s+1}{r+1}\right]  =\left(
0,1\right]  ,
\]
and
\[
\Lambda_{4}=\left[  1,\frac{s+1}{r}\right)  \cap\left[  \frac{s+1}{r+1}
,\frac{s+1}{r}\right)  =\left[  \frac{s+1}{r+1},\frac{s+1}{r}\right)  .
\]
For any $\lambda<\frac{s}{r}$, since $f_{\sigma}$ is monotone increasing in
$\sigma$, the set $f_{\sigma}\left(  \left(  0,1\right]  \right)  $ is
monotone increasing. On the other hand, it is easy to verify that in $\left(
0,\frac{s+1}{r}\right)  $, $f_{0}$ has a unique critical point
\[
1<\lambda_{0}=\frac{\frac{s+1}{r}}{1+\sqrt{\frac{s+1-r}{sr}}}<\frac{s+1}
{r+1},
\]
and since $f_{\sigma}$ is monotone increasing in $\sigma$ for any
$\lambda<\frac{s}{r}$, and monotone decreasing in $\sigma$ for any
$\lambda>\frac{s}{r}$, we have
\[
\max_{\lambda\in\left[  \frac{s+1}{r+1},\frac{s+1}{r}\right)  }f_{\sigma
}\left(  \lambda\right)  =\max_{\lambda\in\left[  \frac{s+1}{r+1},\frac{s}
{r}\right]  }f_{\sigma}\left(  \lambda\right)
\]
which increases with $\sigma$. Hence, $k_{\sigma}$ is monotone increasing.
Next,
\[
\lim_{\sigma\rightarrow0}k_{\sigma}=k_{0}=1
\]
follows from
\[
f_{0}\left(  \frac{s}{r}\right)  =f_{0}\left(  1\right)  =1.
\]
And
\[
\lim_{\sigma\rightarrow\infty}k_{\sigma}=\frac{s+1}{r+1}
\]
follows from
\[
\lim_{\sigma\rightarrow\infty}f_{\sigma}\left(  \lambda\right)  =s+1-\lambda
r.
\]

(iv). When $s<r$, we have
\[
0<\frac{s}{r}<\frac{s+1}{r+1}<1<\frac{s+1}{r},
\]
hence
\[
\Lambda_{3}=\left(  0,\frac{s+1}{r+1}\right]  ,\Lambda_{4}=\left[
1,\frac{s+1}{r}\right)  .
\]
Now the critical point $\lambda_{0}$ of $f_{0}$ satisfies
\[
\frac{s+1}{r+1}<\lambda_{0}<1,
\]
and we have
\[
f_{0}\left(  \frac{s}{r}\right)  =f_{0}\left(  1\right)  =1.
\]
For any $\lambda\geq1>\frac{s}{r}$,
\[
f_{\sigma}\left(  \lambda\right)  \leq f_{0}\left(  \lambda\right)  \leq
f_{0}\left(  1\right)  =f_{0}\left(  \frac{s}{r}\right)  =f_{\sigma}\left(
\frac{s}{r}\right)  ,
\]
hence
\[
f_{\sigma}\left(  \Lambda_{3}\right)  \cap f_{\sigma}\left(  \Lambda
_{4}\right)  =\left(  0,f_{\sigma}\left(  1\right)  \right]  ,
\]
so we have $k_{\sigma}=f_{\sigma}\left(  1\right)  $ which is monotone
decreasing. Limits of $k_{\sigma}$ follow similarly as the case when $s>r$.
\end{proof}

When $d_{2}=d_{1}$, the uniqueness of solutions to $\left(
\ref{Equation stationary}\right)  $ with common sources can be proved directly.

\begin{theorem}
Let $p=r$, $q=s$ and $\left(  \ref{Equation pqrs}\right)  $ be satisfied. Then
$u\equiv u^{\ast},v\equiv v^{\ast}$ is the only solution when $d_{2}=d_{1}$.
\end{theorem}

\begin{proof}
Let $w=u-v$, then we have
\[
\left\{
\begin{array}
[c]{lll}%
d_{1}\triangle w-w+\sigma=0 & \text{in} & \Omega,\\
\frac{\partial w}{\partial\nu}=0 & \text{on} & \partial\Omega,
\end{array}
\right.
\]
hence $w=\sigma$, and $u=v+\sigma$, so we have
\[
\left\{
\begin{array}
[c]{lll}%
d_{1}\triangle v-v+\frac{\left(  v+\sigma\right)  ^{r}}{v^{s}}=0 & \text{in} &
\Omega,\\
\frac{\partial v}{\partial\nu}=0 & \text{on} & \partial\Omega.
\end{array}
\right.
\]
Now maximum principle implies
\[
\bar{v}^{\frac{s+1}{r}}\leq\bar{v}+\sigma\text{ and }\underline{v}^{\frac
{s+1}{r}}\geq\underline{v}+\sigma.
\]
Since
\[
\frac{s+1}{r}>1,
\]
the function
\[
g\left(  t\right)  =t^{\frac{s+1}{r}}-t-\sigma
\]
satisfies
\[
g\left(  0\right)  =-\sigma,\quad g\left(  v^{\ast}\right)  =0,\quad
\lim_{t\rightarrow\infty}g\left(  t\right)  =\infty
\]
and $g$ has only one critical point $t_{1}$ in $\left(  0,\infty\right)  $
which is a local minimum with $g\left(  t_{1}\right)  <-\sigma$. Now $g\left(
\underline{v}\right)  \geq0$ implies $\underline{v}\geq v^{\ast}$ where
$v^{\ast}$ is the unique zero of $g$, and $g\left(  \bar{v}\right)  \leq0$
implies $\bar{v}\leq v^{\ast}$, hence we can conclude $\underline{v}=\bar
{v}=v^{\ast}$.
\end{proof}

\begin{remark}
When $s<r$ and $\sigma>0$, we have $k_{\sigma}<1$, hence the above theorem is
not covered by Corollary \ref{Corollary optimal when sigma>0}.
\end{remark}

\section{Energy method I\label{Section energy 1}}

In this section, we will establish \textit{a priori} bounds independent of
$d_{2}$. The main idea is to use the bound
\[
\int_{\Omega}\frac{u^{r}}{v^{s+1}}\leq\left\vert \Omega\right\vert
\]
to control the nonlinear term in the equation for $u$. This approach can be
traced back to \cite{1986Ni_Takagi} by W. Ni and I. Takagi and
\cite{1987Masuda_Takahashi} by K. Masuda and K. Takahashi. The new ingredient
here is the Moser iteration technique, which enables us to obtain more general result.

We first consider the case when $\sigma=0$.

\begin{theorem}
\label{Theorem upper bounds 0 using MT method} Assume that $\sigma=0$ and
\begin{equation}
\frac{q}{s+1}<\min\left\{  1,\frac{2}{n}\right\}  .\label{Equation
assumption MT00}%
\end{equation}
Then
\[
u\leq c\left(  1+d_{1}^{-\gamma}\right)  ,\quad v\leq c\left(  1+d_{1}
^{-\frac{r}{s+1}\gamma}\right)  ,
\]
where $c,\lambda$ are constants independent of $d_{1},d_{2}$.
\end{theorem}

\begin{proof}
For any $l>1$, multiplying $\left(  \ref{Equation stationary}\right)  _{1}$
with $u^{l-1}$ and integrating over $\Omega$, we obtain
\[
\left(  l-1\right)  \int_{\Omega}d_{1}\left\vert \triangledown u\right\vert
^{2}u^{l-2}+\int_{\Omega}u^{l}=\int_{\Omega}\frac{u^{p+l-1}}{v^{q}}.
\]
Since $\frac{q}{s+1}<1$, applying H\"{o}lder's inequality, we have
\begin{align*}
&  \int_{\Omega}\frac{u^{p+l-1}}{v^{q}}=\int_{\Omega}\left(  \frac{u^{r}
}{v^{s+1}}\right)  ^{\frac{q}{s+1}}\left(  u^{\left(  l+p-1-\frac{qr}
{s+1}\right)  \frac{s+1}{s+1-q}}\right)  ^{1-\frac{q}{s+1}}\\
\leq &  \left(  \int_{\Omega}\frac{u^{r}}{v^{s+1}}\right)  ^{\frac{q}{s+1}
}\left(  \int_{\Omega}u^{\left(  l-\frac{qr-\left(  p-1\right)  \left(
s+1\right)  }{s+1}\right)  \frac{s+1}{s+1-q}}\right)  ^{1-\frac{q}{s+1}}\\
\leq &  \left\vert \Omega\right\vert ^{\frac{q}{s+1}}\left(  \int_{\Omega
}u^{\frac{s+1}{s+1-q}l-\frac{qr-\left(  p-1\right)  \left(  s+1\right)
}{s+1-q}}\right)  ^{1-\frac{q}{s+1}}.
\end{align*}
On the other hand, for any $l\geq\max\left\{  2,4d_{1}\right\}  $, we have
\[
\frac{2d_{1}}{l}\leq d_{1}\frac{4\left(  l-1\right)  }{l^{2}}\leq1.
\]
Lemma \ref{Lemma Sobolev inequality} implies that for any $k\in\left(
1,\frac{2n}{n-2}\right]  $ if $n\geq3$ and $k\in\left(  1,\infty\right)  $ if
$n=1$ or $2$,
\begin{align*}
&  \left(  l-1\right)  \int_{\Omega}d_{1}\left\vert \triangledown u\right\vert
^{2}u^{l-2}+\int_{\Omega}u^{l}\\
=  &  \int_{\Omega}d_{1}\frac{4\left(  l-1\right)  }{l^{2}}\left\vert
\triangledown u^{\frac{l}{2}}\right\vert ^{2}+\int_{\Omega}\left\vert
u^{\frac{l}{2}}\right\vert ^{2}\\
\geq &  c\left(  d_{1}\frac{4\left(  l-1\right)  }{l^{2}}\right)  ^{\frac
{n}{2}-\frac{n}{k}}\left(  \int_{\Omega}\left(  u^{\frac{l}{2}}\right)
^{k}\right)  ^{\frac{2}{k}}\\
\geq &  c\left(  \frac{d_{1}}{l}\right)  ^{\frac{n}{2}-\frac{n}{k}}\left(
\int_{\Omega}\left(  u^{\frac{l}{2}}\right)  ^{k}\right)  ^{\frac{2}{k}}%
\end{align*}
here $c$ is a positive constant depending only on $n$ and $\Omega$. Hence,
\[
c\left(  \frac{d_{1}}{l}\right)  ^{\frac{n}{2}-\frac{n}{k}}\left(
\int_{\Omega}\left(  u^{\frac{l}{2}}\right)  ^{k}\right)  ^{\frac{2}{k}}
\leq\left\vert \Omega\right\vert ^{\frac{q}{s+1}}\left(  \int_{\Omega}
u^{\frac{s+1}{s+1-q}l-\frac{qr-\left(  p-1\right)  \left(  s+1\right)
}{s+1-q}}\right)  ^{1-\frac{q}{s+1}},
\]
i.e.,
\[
\int_{\Omega}u^{\frac{kl}{2}}\leq c\left(  \frac{l}{d_{1}}\right)  ^{\frac
{kn}{4}-\frac{n}{2}}\left(  \int_{\Omega}u^{\frac{s+1}{s+1-q}l-\frac
{qr-\left(  p-1\right)  \left(  s+1\right)  }{s+1-q}}\right)  ^{\frac{k\left(
s+1-q\right)  }{2\left(  s+1\right)  }}.
\]
Setting $l_{0}=2+4d_{1}$, we define $l_{i}$ recursively by
\[
\frac{l_{i}}{2}k=\frac{s+1}{s+1-q}l_{i+1}-\frac{qr-\left(  p-1\right)  \left(
s+1\right)  }{s+1-q},
\]
then
\[
\int_{\Omega}u^{\frac{kl_{i+1}}{2}}\leq c\left(  \frac{l_{i+1}}{d_{1}}\right)
^{\frac{kn}{4}-\frac{n}{2}}\left(  \int_{\Omega}u^{\frac{kl_{i}}{2}}\right)
^{\frac{k\left(  s+1-q\right)  }{2\left(  s+1\right)  }}.
\]
Let
\[
a=\frac{k}{2}\frac{s+1-q}{s+1},\quad b=\frac{qr-\left(  p-1\right)  \left(
s+1\right)  }{s+1},
\]
then
\[
l_{i+1}=al_{i}+b.
\]
Since $\frac{q}{s+1}<\frac{2}{n}$, we can choose $k$ so that $a>1$, and direct
calculation yields
\[
l_{i}=a^{i}\left(  l_{0}+\frac{b}{a-1}\right)  -\frac{b}{a-1}.
\]
Hence
\begin{align*}
&  \int u^{\frac{l_{i+1}}{2}k}\leq\frac{c}{d_{1}^{\frac{kn}{4}-\frac{n}{2}}%
}l_{i+1}^{\frac{kn}{4}-\frac{n}{2}}\left(  \int
u^{\frac{l_{i}}{2}k}\right)
^{a}\\
\leq &  \left(  \frac{c}{d_{1}^{\frac{kn}{4}-\frac{n}{2}}}\right)
^{\sum_{j=0}^{i}a^{j}}%
{\displaystyle\prod\limits_{j=1}^{i+1}} \left(
l_{j}^{\frac{kn}{4}-\frac{n}{2}}\right)  ^{a^{i+1-j}}\left(  \int
u^{\frac{l_{0}}{2}k}\right)  ^{a^{i+1}}\\
\leq &  \left(  \frac{c}{d_{1}^{\frac{kn}{4}-\frac{n}{2}}}\right)
^{\frac{a^{i+1}-1}{a-1}}a^{\left(  \frac{kn}{4}-\frac{n}{2}\right)
\sum
_{j=1}^{i+1}ja^{i+1-j}}\\
&  \left(  l_{0}+\frac{b}{a-1}\right)  ^{\left(  \frac{kn}{4}-\frac{n}%
{2}\right)  \sum_{j=1}^{i+1}a^{i+1-j}}\left(  \int u^{\frac{l_{0}}{2}%
k}\right)  ^{a^{i+1}}\\
=  &  \left(  c\left(  \frac{l_{0}+\frac{b}{a-1}}{d_{1}}\right)
^{\left( \frac{kn}{4}-\frac{n}{2}\right)  }\right)
^{\frac{a^{i+1}-1}{a-1}}a^{\left( \frac{kn}{4}-\frac{n}{2}\right)
\frac{a^{i+2}-\left(  i+2\right) a+i+1}{\left(  a-1\right)
^{2}}}\left(  \int u^{\frac{l_{0}}{2}k}\right) ^{a^{i+1}},
\end{align*}
and
\[
\left(  \int u^{\frac{l_{i+1}}{2}k}\right)  ^{\frac{2}{kl_{i+1}}}\leq\left[
c_{i}\left(  \int u^{\frac{l_{0}}{2}k}\right)  ^{a^{i+1}}\right]  ^{\frac
{2}{k\left(  a^{i+1}\left(  l_{0}+\frac{b}{a-1}\right)  -\frac{b}{a-1}\right)
}}
\]
where
\[
c_{i}=\left(  c\left(  \frac{l_{0}+\frac{b}{a-1}}{d_{1}}\right)  ^{\left(
\frac{kn}{4}-\frac{n}{2}\right)  }\right)  ^{\frac{a^{i+1}-1}{a-1}}a^{\left(
\frac{kn}{4}-\frac{n}{2}\right)  \frac{a^{i+2}-\left(  i+2\right)
a+i+1}{\left(  a-1\right)  ^{2}}}.
\]
Letting $i\rightarrow\infty$, we have
\begin{align*}
&  \left\Vert u\right\Vert _{L^{\infty}\left(  \Omega\right)  }\\
\leq &  \left[  \left(  c\left(  \frac{l_{0}+\frac{b}{a-1}}{d_{1}}\right)
^{\left(  \frac{kn}{4}-\frac{n}{2}\right)  }\right)  ^{\frac{1}{a-1}
}a^{\left(  \frac{kn}{4}-\frac{n}{2}\right)  \frac{a}{\left(  a-1\right)
^{2}}}\int u^{\frac{l_{0}}{2}k}\right]  ^{\frac{2}{k\left(  l_{0}+\frac
{b}{a-1}\right)  }}\\
\leq &  \left[  \left(  c\left(  \frac{l_{0}+\frac{b}{a-1}}{d_{1}}\right)
^{\left(  \frac{kn}{4}-\frac{n}{2}\right)  }\right)  ^{\frac{1}{a-1}
}a^{\left(  \frac{kn}{4}-\frac{n}{2}\right)  \frac{a}{\left(  a-1\right)
^{2}}}\right]  ^{\frac{2}{k\left(  l_{0}+\frac{b}{a-1}\right)  }}\left\Vert
u\right\Vert _{L^{\infty}\left(  \Omega\right)  }^{\frac{l_{0}}{\left(
l_{0}+\frac{b}{a-1}\right)  }}.
\end{align*}
Hence,
\begin{align*}
&  \left\Vert u\right\Vert _{L^{\infty}\left(  \Omega\right)  }\leq\left[
\left(  c\left(  \frac{l_{0}+\frac{b}{a-1}}{d_{1}}\right)  ^{\left(  \frac
{kn}{4}-\frac{n}{2}\right)  }\right)  ^{\frac{1}{a-1}}a^{\left(  \frac{kn}
{4}-\frac{n}{2}\right)  \frac{a}{\left(  a-1\right)  ^{2}}}\right]  ^{\frac
{2}{k}\frac{a-1}{b}}\\
\leq &  c\left(  \frac{l_{0}+\frac{b}{a-1}}{d_{1}}\right)  ^{\left(  \frac
{n}{2}-\frac{n}{k}\right)  \frac{1}{b}}\leq c\left(  1+d_{1}^{-\left(
\frac{n}{2}-\frac{n}{k}\right)  \frac{1}{b}}\right)  .
\end{align*}
The upper bound for $v$ follows from $\left(
\ref{Equation simple inequality 1}\right)  $.
\end{proof}

If $p\geq\frac{qr}{s+1}$, we can actually obtain H\"{o}lder estimate for $u$
directly without using the Moser iteration technique.

\begin{theorem}
Assume in addition to the assumptions of Theorem
\ref{Theorem upper bounds 0 using MT method},
\[
p\geq\frac{qr}{s+1}.
\]
Then for some $\theta\in\left(  0,1\right)  $,
\[
\left\Vert u\right\Vert _{C^{\theta}\left(  \Omega\right)  }\leq c\left(
1+d_{1}^{-\gamma}\right)
\]
where $c,\gamma$ are positive constants independent of $d_{1},d_{2}$.
\end{theorem}

\begin{proof}
We have from the proof of Theorem \ref{Theorem upper bounds 0 using MT method}
, for any $k\in\left(  1,\frac{2n}{n-2}\right]  $ if $n\geq3$, $k\in\left(
1,\infty\right)  $ if $n=1$ or $2$, and for any $l\geq\max\left\{
2,4d_{1}\right\}  $,
\[
\int_{\Omega}u^{\frac{kl}{2}}\leq c\left(  \frac{l}{d_{1}}\right)  ^{\frac
{kn}{4}-\frac{n}{2}}\left(  \int_{\Omega}u^{\frac{s+1}{s+1-q}l-\frac
{qr-\left(  p-1\right)  \left(  s+1\right)  }{s+1-q}}\right)  ^{\frac{k\left(
s+1-q\right)  }{2\left(  s+1\right)  }}.
\]
Setting
\[
k=2\frac{s+1}{s+1-q}
\]
and applying H\"{o}lder's inequality, we have
\begin{align*}
&  \int_{\Omega}u^{\frac{s+1}{s+1-q}l}\leq c\left(  \frac{l}{d_{1}}\right)
^{\frac{nq}{2\left(  s+1-q\right)  }}\int_{\Omega}u^{\frac{s+1}{s+1-q}
l-\frac{qr-\left(  p-1\right)  \left(  s+1\right)  }{s+1-q}}\\
\leq &  c\left(  \frac{l}{d_{1}}\right)  ^{\frac{nq}{2\left(  s+1-q\right)  }
}\left(  \int_{\Omega}u^{\frac{s+1}{s+1-q}l}\right)  ^{1-\frac{qr-\left(
p-1\right)  \left(  s+1\right)  }{\left(  s+1\right)  l}},
\end{align*}
hence
\begin{equation}
\left(  \int_{\Omega}u^{\frac{\left(  s+1\right)  l}{s+1-q}}\right)
^{\frac{s+1-q}{\left(  s+1\right)  l}}\leq c\left(  \frac{l}{d_{1}}\right)
^{\frac{nq}{2\left(  qr-\left(  p-1\right)  \left(  s+1\right)  \right)  }%
}.\label{Equation L^l}%
\end{equation}
Next, for any
\[
\frac{n}{2}<\beta<\frac{s+1}{q},
\]
applying H\"{o}lder's inequality, we have
\begin{align*}
&  \int_{\Omega}\left(  \frac{u^{p}}{v^{q}}\right)  ^{\beta}\leq\left(
\int_{\Omega}\frac{u^{r}}{v^{s+1}}\right)  ^{\beta\frac{q}{s+1}}\left(
\int_{\Omega}u^{\frac{\beta\left(  p-\frac{qr}{s+1}\right)  }{1-\beta\frac
{q}{s+1}}}\right)  ^{1-\beta\frac{q}{s+1}}\\
\leq &  c\left(  1+d_{1}^{-\frac{\beta\left(  p-\frac{qr}{s+1}\right)
nq}{2\left(  qr-\left(  p-1\right)  \left(  s+1\right)  \right)  }}\right)  ,
\end{align*}
where in the last inequality, we used the fact that $u\in L^{\frac{\left(
s+1\right)  l}{s+1-q}}$ for any $l$ large, and
\[
\frac{\beta\left(  p-\frac{qr}{s+1}\right)  }{1-\beta\frac{q}{s+1}}>0.
\]
Now standard elliptic regularity theory yields
\[
\left\Vert u\right\Vert _{W^{2,\beta}\left(  \Omega\right)  }\leq c\left(
\frac{1}{d_{1}}\left\Vert u-\frac{u^{p}}{v^{q}}\right\Vert _{L^{\beta}\left(
\Omega\right)  }+\left\Vert u\right\Vert _{L^{2}\left(  \Omega\right)
}\right)  \leq c\left(  1+d_{1}^{-\gamma}\right)  .
\]
Since $\beta>\frac{n}{2}$, bound for H\"{o}lder norm of $u$ follows from the
Sobolev imbedding.
\end{proof}

\begin{remark}
When $n\geq2$ and $\sigma=0$, W. Ni and I. Takagi \cite{1986Ni_Takagi}
obtained H\"{o}lder bound for $u$ when
\[
\frac{p}{q}=\frac{r}{s+1},\quad\frac{r}{p}>\frac{n}{2}
\]
which is a special case of our theorem.
\end{remark}

Next, we consider the case $\sigma>0$. As in Theorem
\ref{Theorem upper bounds positive copy(1)}, we can relax the conditions on
$p,q,r,s$ with the help of trivial lower bounds.

\begin{theorem}
\label{Theorem NiTakagi positive}Let $\sigma>0$ and
\[
\frac{p-1}{r}<\min\left\{  1,\frac{2}{n}\right\}  .
\]
Then we have
\[
u\leq c\left(  1+d_{1}^{-\gamma}\right)  ,\quad v\leq c\left(  1+d_{1}
^{-\frac{r}{s+1}\gamma}\right)  ,
\]
where $c,\gamma$ are positive constants independent of $d_{1},d_{2}$.
\end{theorem}

\begin{proof}
We choose $\varepsilon\geq0$, such that
\[
\frac{p-1}{r}<\frac{q-\varepsilon}{s+1}<\min\left\{  1,\frac{2}{n}\right\}  .
\]
For any $l>1$, we have the energy estimate
\[
\left(  l-1\right)  \int_{\Omega}d_{1}\left\vert \triangledown u\right\vert
^{2}u^{l-2}+\int_{\Omega}u^{l}\leq\int_{\Omega}\frac{u^{p+l-1}}{v^{q}}
+\sigma\int_{\Omega}u^{l-1}.
\]
If
\[
\left(  l-1\right)  \int_{\Omega}d_{1}\left\vert \triangledown u\right\vert
^{2}u^{l-2}+\int_{\Omega}u^{l}\leq2\sigma\int_{\Omega}u^{l-1}\leq
2\sigma\left\vert \Omega\right\vert ^{\frac{1}{l}}\left(  \int_{\Omega}
u^{l}\right)  ^{\frac{l-1}{l}},
\]
then
\[
\int_{\Omega}u^{l}\leq2\sigma\int_{\Omega}u^{l-1}\leq2\sigma\left\vert
\Omega\right\vert ^{\frac{1}{l}}\left(  \int_{\Omega}u^{l}\right)
^{\frac{l-1}{l}},
\]
hence
\[
\int_{\Omega}u^{l}\leq\left(  2\sigma\right)  ^{l}\left\vert \Omega\right\vert
.
\]
Otherwise,
\[
\left(  l-1\right)  \int_{\Omega}d_{1}\left\vert \triangledown u\right\vert
^{2}u^{l-2}+\int_{\Omega}u^{l}\leq2\int_{\Omega}\frac{u^{p+l-1}}{v^{q}}
\leq2\sigma^{\frac{r}{s+1}\varepsilon}\int_{\Omega}\frac{u^{p+l-1}
}{v^{q-\varepsilon}},
\]
the arguments leading to $\left(  \ref{Equation L^l}\right)  $ implies that
for any $l\geq\max\left\{  2,4d_{1}\right\}  $,
\[
\left(  \int_{\Omega}u^{\frac{\left(  s+1\right)  l}{s+1-\left(
q-\varepsilon\right)  }}\right)  ^{\frac{s+1-\left(  q-\varepsilon\right)
}{\left(  s+1\right)  l}}\leq c\left(  \frac{l}{d_{1}}\right)  ^{\frac
{n\left(  q-\varepsilon\right)  }{2\left(  \left(  q-\varepsilon\right)
r-\left(  p-1\right)  \left(  s+1\right)  \right)  }},
\]
here $c$ is a constant independent of $d_{1},d_{2}$. Hence, in any case, we
have for $l$ large
\[
\left(  \int_{\Omega}u^{l}\right)  ^{\frac{1}{l}}\leq c\left(  1+d_{1}
^{-\gamma}\right)
\]
where $c,\gamma$ are positive constants independent of $d_{1},d_{2}$. Next,
\[
\frac{u^{p}}{v^{q}}\leq\sigma^{-q}u^{p},
\]
hence $\frac{u^{p}}{v^{q}}\in L^{\beta}$ for any $\beta$. Choosing
$\beta>\frac{n}{2}$, elliptic regularity theory yields
\[
u\in W^{2,\beta}\hookrightarrow C^{\theta}\left(  \overline{\Omega}\right)  ,
\]
especially, we have $u\leq c\left(  1+d_{1}^{-\lambda}\right)  $. Upper bound
for $v$ follows from $\left(  \ref{Equation simple inequality 1}\right)  $.
\end{proof}

\begin{proof}
[Proof of Theorem \ref{Theorem GM estimate}]Part (i) is exactly Theorem
\ref{Theorem upper bounds 0 using MT method}; Part (ii)
\[
\frac{p-1}{r}<\frac{2}{n}
\]
follows from Theorem \ref{Theorem NiTakagi positive}.
\end{proof}

\section{Energy Method II\label{Section energy 2}}

In this section, we will modify the approach in \cite{1994DelPino} by M. del
Pino to establish \textit{a priori} bounds which are uniform when $d_{1}
,d_{2}$ are large. The main idea is to use $L^{1}$ norm to control upper and
lower bounds for $u,v$.

We first look at the case when $\sigma=0$.

\begin{theorem}
\label{Theorem a family of lower bounds}Assume that $\sigma=0$, $r<\frac
{n}{n-2}$ and there exists $\delta\in\left(  0,1\right]  $ such that
\begin{align*}
0  &  <\frac{1-\delta}{r}+\frac{\delta}{p}<1,\\
\frac{\frac{\left(  1-\delta\right)  s}{r}+\frac{\delta q}{p}}{\frac
{r-1+\delta}{r}-\frac{\delta}{p}}  &  <\frac{n}{n-2}\text{ or }\frac
{\frac{\left(  1-\delta\right)  s}{r}+\frac{\delta q}{p}}{\frac{r-1+\delta}
{r}-\frac{\delta}{p}}\leq s+1.
\end{align*}
Then for any $\eta>0$, there exists positive constant $c=c\left(
p,q,r,s,\eta\right)  $ such that
\[
u\geq c,\quad v>c^{\frac{s+1}{r}}
\]
whenever $d_{1},d_{2}\geq\eta$.
\end{theorem}

\begin{proof}
Under the condition
\[
0<\delta\leq1,\quad0<\frac{1-\delta}{r}+\frac{\delta}{p}<1,
\]
applying H\"{o}lder's inequality, we have
\[
\int_{\Omega}u\leq\left(  \int_{\Omega}\frac{u^{r}}{v^{s}}\right)
^{\frac{1-\delta}{r}}\left(  \int_{\Omega}\frac{u^{p}}{v^{q}}\right)
^{\frac{\delta}{p}}\left(  \int_{\Omega}v^{\frac{\left(  \frac{1-\delta}
{r}\right)  s+\frac{\delta q}{p}}{\frac{r-1+\delta}{r}-\frac{\delta}{p}}
}\right)  ^{\frac{r-1+\delta}{r}-\frac{\delta}{p}}.
\]
Here if $\delta=1$, we simply drop the term on $\int_{\Omega}\frac{u^{r}
}{v^{s}}$. Using Lemmas \ref{Lemma trivial estimates} and
\ref{Lemma reverse Holder}, we have
\[
\int_{\Omega}\frac{u^{r}}{v^{s}}=\int_{\Omega}v\leq c\left(  \int_{\Omega
}v^{s+1}\right)  ^{\frac{1}{s+1}}\leq c\left(  \int_{\Omega}u^{r}\right)
^{\frac{1}{s+1}}\leq c\left(  \int_{\Omega}u\right)  ^{\frac{r}{s+1}},
\]
where we have used the assumption $r<\frac{n}{n-2}$.

Next, if
\[
\frac{\left(  \frac{1-\delta}{r}\right)  s+\frac{\delta q}{p}}{\frac
{r-1+\delta}{r}-\frac{\delta}{p}}<\frac{n}{n-2},
\]
we can apply Lemma \ref{Lemma reverse Holder}, and obtain
\[
\int_{\Omega}v^{\frac{\left(  \frac{1-\delta}{r}\right)  s+\frac{\delta q}{p}
}{\frac{r-1+\delta}{r}-\frac{\delta}{p}}}\leq c\left(  \int_{\Omega}v\right)
^{\frac{\left(  \frac{1-\delta}{r}\right)  s+\frac{\delta q}{p}}
{\frac{r-1+\delta}{r}-\frac{\delta}{p}}}\leq c\left(  \int_{\Omega}u\right)
^{\frac{r}{s+1}\frac{\left(  \frac{1-\delta}{r}\right)  s+\frac{\delta q}{p}
}{\frac{r-1+\delta}{r}-\frac{\delta}{p}}},
\]
and if
\[
\frac{\left(  \frac{1-\delta}{r}\right)  s+\frac{\delta q}{p}}{\frac
{r-1+\delta}{r}-\frac{\delta}{p}}\leq s+1,
\]
we can apply H\"{o}lder's inequality,
\[
\int_{\Omega}v^{\frac{\left(  \frac{1-\delta}{r}\right)  s+\frac{\delta q}{p}
}{\frac{r-1+\delta}{r}-\frac{\delta}{p}}}\leq c\left(  \int_{\Omega}
v^{s+1}\right)  ^{\frac{1}{s+1}\frac{\left(  \frac{1-\delta}{r}\right)
s+\frac{\delta q}{p}}{\frac{r-1+\delta}{r}-\frac{\delta}{p}}}\leq c\left(
\int_{\Omega}u\right)  ^{\frac{r}{s+1}\frac{\left(  \frac{1-\delta}{r}\right)
s+\frac{\delta q}{p}}{\frac{r-1+\delta}{r}-\frac{\delta}{p}}}.
\]
Finally,
\[
\int_{\Omega}\frac{u^{p}}{v^{q}}=\int_{\Omega}u.
\]
Combining all these estimates, we have
\[
\int_{\Omega}u\leq c\left(  \int_{\Omega}u\right)  ^{1+\frac{\delta}{p\left(
s+1\right)  }\left(  qr-\left(  p-1\right)  \left(  s+1\right)  \right)  }.
\]
Since
\[
\frac{\delta}{p\left(  s+1\right)  }\left(  qr-\left(  p-1\right)  \left(
s+1\right)  \right)  >0,
\]
we have
\[
\int_{\Omega}u\geq c.
\]
Now Lemma \ref{Lemma u,v larger than their average} implies,
\[
\underline{u}=\min_{x\in\overline{\Omega}}u\geq c,
\]
and from $\left(  \ref{Equation simple inequality 1}\right)  $,
\[
\underline{v}\geq\underline{u}^{\frac{r}{s+1}}\geq c^{\frac{r}{s+1}}.
\]

\end{proof}

\begin{remark}
The assumptions in theorem \ref{Theorem a family of lower bounds} seem
complicated. However, these assumptions automatically hold when the dimension
$n=2$.
\end{remark}

Since $\delta$ is a free parameter, a family of \textit{a priori} lower bounds
can be deduced from Theorem \ref{Theorem a family of lower bounds}.

\begin{corollary}
Assume $\sigma=0$, $r<\frac{n}{n-2}$ and
\[
\frac{q}{p-1}<\frac{n}{n-2}\text{ or }\frac{q}{p-1}\leq s+1.
\]
Given $\eta>0$, there exists a positive constant $c=c\left(  p,q,r,s,\eta
\right)  $ such that
\[
u\geq c,\quad v>c^{\frac{s+1}{r}}
\]
whenever $d_{1},d_{2}\geq\eta$.
\end{corollary}

\begin{proof}
Let $\delta=1$, it is straightforward that our conditions are equivalent to
the assumptions in Theorem \ref{Theorem a family of lower bounds}.
\end{proof}

\begin{corollary}
Assume $\sigma=0$, $1<r<\frac{n}{n-2}$ and
\[
\frac{s}{r-1}<\frac{n}{n-2}\text{ or }\frac{s}{r-1}<s+1.
\]
Given $\eta>0$, there exists a positive constant $c=c\left(  p,q,r,s,\eta
\right)  $ such that
\[
u\geq c,\quad v\geq c^{\frac{s+1}{r}}
\]
whenever $d_{1},d_{2}\geq\eta$.
\end{corollary}

\begin{proof}
Under our assumptions, it is easy to check that the assumptions in Theorem
\ref{Theorem a family of lower bounds} will be satisfied for sufficiently
small $\delta>0$.
\end{proof}

\begin{remark}
In \cite{1994DelPino}, it was proved that $1<r<\frac{n}{n-2}$ and $\frac
{s}{r-1}<\frac{n}{n-2}$ imply positive lower a priori bounds which is a
special case of the above corollary.
\end{remark}

Next, we consider the \textit{a priori} upper bounds of $u,v$. We have two
different methods to get such bounds which were discussed respectively in
Section \ref{Section maximum principle} and Section \ref{Section energy 1}.

\begin{theorem}
\label{Theorem 0 energy 2}Let $\sigma=0$,
\[
\frac{p-1}{r}<\min\left(  1,\frac{2}{n}\right)
\]
and assume that $v\geq c_{1}$. Then
\[
u\leq c_{2},\quad v\leq c_{2}^{\frac{r}{s+1}}
\]
where $c_{2}$ is a positive constant depending on $n,p,q,r,s,d_{1}$ and
$c_{1}$.
\end{theorem}

\begin{proof}
See the proof of Theorem \ref{Theorem NiTakagi positive}.
\end{proof}

\begin{proof}
[Proof of Theorem \ref{Theorem GM estimate 1}]Part (i) follows from
\ref{Theorem a family of lower bounds}; Part (ii) follows from
\ref{Theorem 0 energy 2}.
\end{proof}

When $n=2$, our theorems imply

\begin{corollary}
\label{Corollary n=2}Let $n=2$ and $\sigma=0$. Given $\eta>0$, there exists a
positive constant $c_{1}$ depending only on $p,q,r,s,\eta$ and $\Omega$, such
that
\[
u\geq c_{1},\quad v\geq c_{1}^{\frac{r}{s+1}}
\]
whenever $d_{1},d_{2}>\eta$. If in addition
\[
\frac{p-1}{r}<1,
\]
then there exists a positive constant $c_{2}$ depending only on $p,q,r,s,\eta$
and $\Omega$, such that
\[
u\leq c_{2},\quad v\leq c_{2}^{\frac{r}{s+1}}
\]
whenever $d_{1},d_{2}>\eta$.
\end{corollary}

\section{Nonexistence of nontrivial solutions via energy
estimate\label{Section nonexistence}}

It is generally expected that uniform bounds imply nonexistence of nontrivial
solutions when the diffusion constants are sufficiently large. More precisely,
we have the following result.

\begin{theorem}
\label{Theorem nonexistence when d_1 is large}Let $D\subset\left(
0,\infty\right)  \times\left(  0,\infty\right)  $ be nonempty. Assuming that
we have positive lower and upper a priori bounds for solutions to $\left(
\ref{Equation stationary}\right)  $ which are uniform for $\left(  d_{1}
,d_{2}\right)  \in D$. Then there exists a positive constant $C$ such that for
any $\left(  d_{1},d_{2}\right)  \in D$ satisfying
\[
d_{1}\geq C,
\]
$\left(  u,v\right)  \equiv\left(  u^{\ast},v^{\ast}\right)  $ is the only
solution of $\left(  \ref{Equation stationary}\right)  $.
\end{theorem}

\begin{proof}
Let
\[
\bar{u}=\frac{1}{\left\vert \Omega\right\vert }\int_{\Omega}u\left(  x\right)
dx,\quad\bar{v}=\frac{1}{\left\vert \Omega\right\vert }\int_{\Omega}v\left(
x\right)  dx.
\]
First, multiplying the second equation of $\left(  \ref{Equation stationary}
\right)  $ with $v-\bar{v}$ and integrating over $\Omega$, we have
\begin{align*}
&  d_{2}\int_{\Omega}\left\vert \triangledown v\right\vert ^{2}+\int_{\Omega
}\left\vert v-\bar{v}\right\vert ^{2}=\int_{\Omega}\left(  \frac{u^{r}}{v^{s}
}-\frac{\bar{u}^{r}}{\bar{v}^{s}}\right)  \left(  v-\bar{v}\right) \\
=  &  \int_{\Omega}\left(  \frac{u^{r}}{v^{s}}-\frac{u^{r}}{\bar{v}^{s}
}\right)  \left(  v-\bar{v}\right)  +\int_{\Omega}\left(  \frac{u^{r}}{\bar
{v}^{s}}-\frac{\bar{u}^{r}}{\bar{v}^{s}}\right)  \left(  v-\bar{v}\right) \\
\leq &  \int_{\Omega}\left(  \frac{u^{r}}{\bar{v}^{s}}-\frac{\bar{u}^{r}}
{\bar{v}^{s}}\right)  \left(  v-\bar{v}\right)  \leq C\int_{\Omega}\left(
u-\bar{u}\right)  \left(  v-\bar{v}\right) \\
\leq &  \frac{1}{2}\int_{\Omega}\left\vert v-\bar{v}\right\vert ^{2}
+C\int_{\Omega}\left\vert u-\bar{u}\right\vert ^{2}%
\end{align*}
hence
\begin{equation}
\int_{\Omega}\left\vert v-\bar{v}\right\vert ^{2}\leq C\int_{\Omega}\left\vert
u-\bar{u}\right\vert ^{2}.\label{Equation
constant}%
\end{equation}
Next, from the first equation of $\left(  \ref{Equation stationary}\right)  $,
we have
\begin{align*}
&  d_{1}\int_{\Omega}\left\vert \triangledown u\right\vert ^{2}+\int_{\Omega
}\left\vert u-\bar{u}\right\vert ^{2}=\int_{\Omega}\left(  \frac{u^{p}}{v^{q}
}-\frac{\bar{u}^{p}}{\bar{v}^{q}}\right)  \left(  u-\bar{u}\right) \\
=  &  \int_{\Omega}\left(  \frac{u^{p}}{v^{q}}-\frac{\bar{u}^{p}}{v^{q}
}\right)  \left(  u-\bar{u}\right)  +\left(  \frac{\bar{u}^{p}}{v^{q}}
-\frac{\bar{u}^{p}}{\bar{v}^{q}}\right)  \left(  u-\bar{u}\right) \\
\leq &  C\int_{\Omega}\left(  u-\bar{u}\right)  ^{2}+\int_{\Omega}\left\vert
u-\bar{u}\right\vert \cdot\left\vert v-\bar{v}\right\vert \\
\leq &  C\int_{\Omega}\left(  u-\bar{u}\right)  ^{2}+\int_{\Omega}\left\vert
v-\bar{v}\right\vert ^{2}\leq C\int_{\Omega}\left(  u-\bar{u}\right)  ^{2}.
\end{align*}
Hence
\[
d_{1}\left\Vert \triangledown u\right\Vert _{L^{2}\left(  \Omega\right)  }
^{2}\leq C\int_{\Omega}\left(  u-\bar{u}\right)  ^{2}\leq C\left\Vert
\triangledown u\right\Vert _{L^{2}\left(  \Omega\right)  }^{2}.
\]
If $d_{1}\geq C$, we have
\[
\left\Vert \triangledown u\right\Vert _{L^{2}\left(  \Omega\right)  }=0,
\]
and $u$ is constant. Now $\left(  \ref{Equation constant}\right)  $ implies
that $v$ is also constant, and the only constant solution to $\left(
\ref{Equation stationary}\right)  $ is $\left(  u,v\right)  \equiv\left(
u^{\ast},v^{\ast}\right)  $.
\end{proof}

Now we are ready to prove Theorem \ref{Theorem nonexistence}.

\begin{proof}
[Proof of Theorem \ref{Theorem nonexistence}]Part (i) follows from Theorem
\ref{Theorem upper bounds positive copy(1)}; (ii) follows from Theorems
\ref{Theorem NiTakagi positive}; (iii) follows from Corollary
\ref{Corollary n=2}.
\end{proof}

\section{Existence of nontrivial solutions\label{Section existence}}

In this section, we will use topological degree theory to show the existence
of nontrivial solutions to $\left(  \ref{Equation stationary}\right)  $ under
suitable conditions. This approach has been used by many authors, for example,
M. del Pino in \cite{1994DelPino}.

Let $X=C^{0}\left(  \overline{\Omega}\right)  \times C^{0}\left(
\overline{\Omega}\right)  $ be the Banach space with norm
\[
\left\Vert \left(  u,v\right)  \right\Vert _{X}=\max\left\{  \left\Vert
u\right\Vert _{L^{\infty}\left(  \Omega\right)  },\left\Vert v\right\Vert
_{L^{\infty}\left(  \Omega\right)  }\right\}
\]
and $X^{+}$ be the positive cone in $X$, i.e.,
\[
X^{+}=\left\{  \left(  u,v\right)  \in X:u>0,v>0\text{ in }\overline{\Omega
}\right\}  .
\]
We define solution operators $S=\left(  I-d_{2}\triangle\right)  ^{-1}$ and
$R=\left(  \varpi I-d_{1}\triangle\right)  ^{-1}$ under Neumann boundary
conditions. Here $\varpi>0$ is a large constant to be determined later. Let
\[
T\left(  u,v\right)  =\left(  R\left(  f\left(  u,v\right)  +\varpi u\right)
,S\left(  g\left(  u,v\right)  +v\right)  \right)
\]
where
\[
f\left(  u,v\right)  =-u+\frac{u^{p}}{v^{q}}+\sigma,g\left(  u,v\right)
=-v+\frac{u^{r}}{v^{s}}.
\]
Then $T$ is an operator defined on $X^{+}$ and it is easy to check that
$\left(  u,v\right)  $ is a positive solution to $\left(
\ref{Equation stationary}\right)  $ if and only if it is a fixed point of $T$
in $X^{+}$, i.e.,
\[
T\left(  u,v\right)  =\left(  u,v\right)  .
\]

Now we consider the linearization of $\left(  \ref{Equation stationary}
\right)  $ around $\left(  u^{\ast},v^{\ast}\right)  $,
\begin{equation}
\left\{
\begin{array}
[c]{lll}%
d_{1}\triangle h+f_{u}\left(  u^{\ast},v^{\ast}\right)  h+f_{v}\left(
u^{\ast},v^{\ast}\right)  k=0 & \text{in} & \Omega,\\
d_{2}\triangle k+g_{u}\left(  u^{\ast},v^{\ast}\right)  h+g_{v}\left(
u^{\ast},v^{\ast}\right)  k=0 & \text{in} & \Omega,\\
\frac{\partial h}{\partial\nu}=\frac{\partial k}{\partial\nu}=0 & \text{on} &
\partial\Omega.
\end{array}
\right. \label{Equation linear system}%
\end{equation}
Let $0=\lambda_{0}<\lambda_{1}<\lambda_{2}<\cdots$ be the eigenvalues of
$-\triangle$ under Neumann boundary conditions in $\Omega$. We also use
$m_{i}$ to denote the multiplicity of eigenvalue $\lambda_{i}$,
$i=0,1,2,\cdots$. In the remaining part of the paper, we will simply use
$f_{u}$ to denote $f_{u}\left(  u^{\ast},v^{\ast}\right)  $, the same applies
to $f_{v}$, $g_{u}$ and $g_{v}$. The linear system $\left(
\ref{Equation linear system}\right)  $\ will possess a nontrivial solution if
and only if the matrix
\[
\left(
\begin{array}
[c]{cc}%
f_{u}-d_{1}\lambda_{i} & f_{v}\\
g_{u} & g_{v}-d_{2}\lambda_{i}%
\end{array}
\right)  =\left(
\begin{array}
[c]{cc}%
p\left(  u^{\ast}\right)  ^{p-1-\frac{qr}{s+1}}-1-d_{1}\lambda_{i} & -q\left(
u^{\ast}\right)  ^{p-\frac{\left(  q+1\right)  r}{s+1}}\\
r\left(  u^{\ast}\right)  ^{\frac{r}{s+1}-1} & -s-1-d_{2}\lambda_{i}%
\end{array}
\right)
\]
is singular for some $i$. Hence, given $d_{2}>0$, for each $i\geq1$, the
linear system $\left(  \ref{Equation linear system}\right)  $\ will possess a
nontrivial solution if and only if
\[
d_{1}=d_{1i}\equiv\frac{1}{\lambda_{i}}\left[  p\left(  u^{\ast}\right)
^{p-1-\frac{qr}{s+1}}-1-\frac{qr\left(  u^{\ast}\right)  ^{p-1-\frac{qr}{s+1}
}}{s+1+d_{2}\lambda_{i}}\right]  .
\]
We also define for any $d>0$,
\[
A_{d}=\left\{  i\geq1:d<d_{1i}\right\}  ,N_{d}=\sum_{i\in A_{d}}m_{i}.
\]

\begin{lemma}
If $d_{1}\neq d_{1i}$, $i=1,2,3,\cdots$, then for every sufficiently small
neighborhood $V$ of $\left(  u^{\ast},v^{\ast}\right)  $, $T$ has no fixed
point on $\partial V$ and
\[
\deg\left(  I-T,V,\left(  0,0\right)  \right)  =\left(  -1\right)  ^{N_{d_{1}
}}
\]
provided that $\varpi$ is sufficiently large.
\end{lemma}

\begin{proof}
Let $L$ be the Fr\'{e}chet derivative of $T$ at $\left(  u^{\ast},v^{\ast
}\right)  $. For any $\left(  h,k\right)  \in X$, we have
\[
L\left(  h,k\right)  =\left(  R\left(  f_{u}h+f_{v}k+\varpi h\right)
,S\left(  g_{u}h+g_{v}k+k\right)  \right)  .
\]
If $I-L$ is nonsingular, then $\left(  u^{\ast},v^{\ast}\right)  $ is an
isolated fixed point of $T$ and, for every sufficiently small neighborhood $V$
of $\left(  u^{\ast},v^{\ast}\right)  $,
\[
\deg\left(  I-T,V,\left(  0,0\right)  \right)  =\left(  -1\right)  ^{\eta}
\]
where $\eta$ is the number of negative eigenvalues counting algebraic
multiplicities of $I-L$. So we only need to show $I-L$ is nonsingular and
$\eta=N_{d_{1}}$. Observe that $-\mu\leq0$ is an eigenvalue if and only if the
system
\begin{equation}
\left\{
\begin{array}
[c]{lll}%
-\left(  \mu+1\right)  d_{1}\triangle h+\mu\varpi h=f_{u}h+f_{v}k, & \text{in}
& \Omega,\\
-\left(  \mu+1\right)  d_{2}\triangle k+\mu k=g_{u}h+g_{v}k & \text{in} &
\Omega,\\
\frac{\partial h}{\partial\nu}=\frac{\partial k}{\partial\nu}=0 & \text{on} &
\partial\Omega
\end{array}
\right. \label{Equation linear for mu}%
\end{equation}
has nontrivial solutions, which is equivalent to that the matrix
\begin{align*}
&  \left(
\begin{array}
[c]{cc}%
f_{u}-\varpi\mu-\left(  \mu+1\right)  d_{1}\lambda_{i} & f_{v}\\
g_{u} & g_{v}-\mu-\left(  \mu+1\right)  d_{2}\lambda_{i}%
\end{array}
\right) \\
=  &  \left(
\begin{array}
[c]{cc}%
p\left(  u^{\ast}\right)  ^{p-1-\frac{qr}{s+1}}-1-\varpi\mu-\left(
\mu+1\right)  d_{1}\lambda_{i} & -q\left(  u^{\ast}\right)  ^{p-\frac{\left(
q+1\right)  r}{s+1}}\\
r\left(  u^{\ast}\right)  ^{\frac{r}{s+1}-1} & -s-1-\mu-\left(  \mu+1\right)
d_{2}\lambda_{i}%
\end{array}
\right)
\end{align*}
is singular for some $i\geq0$, i.e.,
\begin{align}
&  \left(  \varpi+d_{1}\lambda_{i}\right)  \left(  1+d_{2}\lambda_{i}\right)
\mu^{2}\nonumber\\
+  &  \left[  \left(  s+1+d_{2}\lambda_{i}\right)  \left(  \varpi+d_{1}
\lambda_{i}\right)  +\left(  1+d_{2}\lambda_{i}\right)  \left(  -p\left(
u^{\ast}\right)  ^{p-1-\frac{qr}{s+1}}+1+d_{1}\lambda_{i}\right)  \right]
\mu\nonumber\\
+  &  qr\left(  u^{\ast}\right)  ^{p-1-\frac{qr}{s+1}}-\left[  p\left(
u^{\ast}\right)  ^{p-1-\frac{qr}{s+1}}-1-d_{1}\lambda_{i}\right]  \left[
s+1+d_{2}\lambda_{i}\right]  =0.\label{Equation mu}%
\end{align}
If we choose $\varpi$ sufficiently large, the left hand side of $\left(
\ref{Equation mu}\right)  $ is monotone increasing in $\mu\geq0$. When $i=0$,
we have $\lambda_{i}=0$, using $u^{\ast}\geq1$ and $\left(
\ref{Equation pqrs}\right)  $, it is easy to check
\[
qr\left(  u^{\ast}\right)  ^{p-1-\frac{qr}{s+1}}-\left(  p\left(  u^{\ast
}\right)  ^{p-1-\frac{qr}{s+1}}-1\right)  \left(  s+1\right)  >0,
\]
hence the matrix can't be singular for any $\mu\geq0$. Now for any $i\geq1$,
we can solve $d_{1}$ from $\left(  \ref{Equation mu}\right)  $,
\[
d_{1}=\frac{p\left(  u^{\ast}\right)  ^{p-1-\frac{qr}{s+1}}-1-\varpi\mu
}{\left(  \mu+1\right)  \lambda_{i}}-\frac{qr\left(  u^{\ast}\right)
^{p-1-\frac{qr}{s+1}}}{\left(  \mu+1\right)  \lambda_{i}\left(  s+1+\mu
+\left(  \mu+1\right)  d_{2}\lambda_{i}\right)  }\equiv p_{i}\left(
\mu\right)  .
\]
When $\varpi$ is sufficiently large, one can check that $p_{i}\left(
\mu\right)  $ is monotone decreasing for any $\mu\geq0$, and
\[
p_{i}\left(  0\right)  =d_{1i},\lim_{\mu\rightarrow\infty}p_{i}\left(
\mu\right)  =-\frac{\varpi}{\lambda_{i}}.
\]
Since $d_{1}\neq d_{1i}$, we conclude that $\mu=0$ is not an eigenvalue to
$\left(  \ref{Equation linear for mu}\right)  $, hence $I-L$ is nonsingular.
If $d_{1}<d_{1i}$, then there exists a unique $\mu>0$ such that
\[
d_{1}=p_{i}\left(  \mu\right)  .
\]
And each eigenfunction of $\lambda_{i}$ gives rise to an eigenfunction of
$\left(  \ref{Equation linear for mu}\right)  $. Hence, we have $\eta
=N_{d_{1}}$.
\end{proof}

Next, we consider a one-parameter family of elliptic systems
\begin{equation}
\left\{
\begin{array}
[c]{lll}%
d_{1}\triangle u-u+\tau\left(  \frac{u^{p}}{v^{q}}+\sigma\right)  +\left(
1-\tau\right)  \rho=0 & \text{in} & \Omega,\\
d_{2}\triangle v-v+\tau\frac{u^{r}}{v^{s}}+\left(  1-\chi_{\tau}\right)
\rho=0 & \text{in} & \Omega,\\
\frac{\partial u}{\partial\nu}=\frac{\partial v}{\partial\nu}=0 & \text{on} &
\partial\Omega
\end{array}
\right. \label{Equation deformation}%
\end{equation}
with parameter $\tau\in\left[  0,1\right]  $. (We have abused the notation
here since the parameter $\tau$ has nothing to do with the response rate in
$\left(  \ref{GM PDE system}\right)  $.) In $\left(
\ref{Equation deformation}\right)  $, $\rho$ is a given positive constant and
\[
\chi_{\tau}=\left\{
\begin{array}
[c]{ccc}%
2\tau & \text{if} & \tau\in\left[  0,\frac{1}{2}\right]  ,\\
1 & \text{if} & \tau\in\left[  \frac{1}{2},1\right]  .
\end{array}
\right.
\]
When $\tau$ changes from $0$ to $1$, $\left(  \ref{Equation deformation}
\right)  $ serves as a deformation from a trivial system which has a unique
solution $\left(  u,v\right)  \equiv\left(  \rho,\rho\right)  $ to $\left(
\ref{Equation stationary}\right)  $. Our choice of $\chi_{\tau}$ simplifies
the proof of \textit{a priori} estimates in Proposition
\ref{Proposition delta}.

For $\eta>0$, we denote
\[
\Lambda_{\eta}=\left\{  \left(  u,v\right)  \in X:\eta<u,v<\frac{1}{\eta
}\text{ in }\overline{\Omega}\right\}  .
\]

\begin{lemma}
\label{Lemma index 1}Assume that positive solutions to $\left(
\ref{Equation deformation}\right)  $ satisfies a priori bound
\begin{equation}
0<\alpha\leq u,v\leq\beta\label{Equation uniform bound}%
\end{equation}
for some positive constants $\alpha,\beta$ independent of $\tau$. Then there
exists $\eta>0$ such that $T$ has no fixed point on $\partial\Lambda_{\eta}$
and
\[
\deg\left(  I-T,\Lambda_{\eta},\left(  0,0\right)  \right)  =1.
\]

\end{lemma}

\begin{proof}
Let
\[
T_{\tau}\left(  u,v\right)  =\left(  R_{\tau}\left(  \tau f\left(  u,v\right)
+\left(  1-\tau\right)  \rho+\tau\varpi u\right)  ,S\left(  \tau\left(
g\left(  u,v\right)  +v\right)  +\left(  1-\chi_{\tau}\right)  \rho\right)
\right)
\]
where
\[
R_{\tau}=\left(  \left(  \tau\varpi+1-\tau\right)  I-d_{1}\triangle\right)
^{-1}.
\]
Then for each $\tau\in\left[  0,1\right]  $, $T_{\tau}$ is a compact operator
from $\Lambda_{\eta}$ into $X$. Furthermore, $T_{1}=T$ and $T_{0}\equiv\left(
\rho,\rho\right)  $. It is easy to check that $\left(  u,v\right)  $ is a
fixed point of $T_{\tau}$ if and only if it is a solution to $\left(
\ref{Equation deformation}\right)  $. The bounds in $\left(
\ref{Equation uniform bound}\right)  $ then implies that for $\eta$
sufficiently small, for each $\tau\in\left[  0,1\right]  $, $T_{\tau}$ has no
fixed point on $\partial\Lambda_{\eta}$. Hence,
\[
\deg\left(  I-T_{1},\Lambda_{\eta},\left(  0,0\right)  \right)  =\deg\left(
I-T_{0},\Lambda_{\eta},\left(  0,0\right)  \right)  ,
\]
i.e.,
\[
\deg\left(  I-T,\Lambda_{\eta},\left(  0,0\right)  \right)  =\deg\left(
I-\left(  \rho,\rho\right)  ,\Lambda_{\eta},\left(  0,0\right)  \right)
=\deg\left(  I,\Lambda_{\eta},\left(  \rho,\rho\right)  \right)  =1
\]
since $\left(  \rho,\rho\right)  \in\Lambda_{\eta}$.
\end{proof}

\begin{theorem}
\label{Theorem existence}Under the assumption of Lemma \ref{Lemma index 1}. If
$d_{1}\neq d_{1i}$, $i=1,2,3,\cdots$, and $N_{d_{1}}$ is odd, then there
exists at least one nontrivial positive solution to $\left(
\ref{Equation stationary}\right)  $.
\end{theorem}

\begin{proof}
From the properties of topological degree, we have
\[
\deg\left(  I-T,\Lambda_{\eta}\backslash\overline{V},\left(  0,0\right)
\right)  =1-\left(  -1\right)  ^{N_{d_{1}}}=2\neq0,
\]
hence $T$ has at least one fixed point in $\Lambda_{\eta}\backslash
\overline{V}$, which is a nontrivial solution to $\left(
\ref{Equation stationary}\right)  $.
\end{proof}

\begin{remark}
One necessary condition to apply Theorem \ref{Theorem existence} is that
$\sigma$ is small. Actually, if
\begin{equation}
\sigma\geq\left(  p-1\right)  p^{\frac{s+1}{qr-\left(  p-1\right)  \left(
s+1\right)  }-1},\label{Equation smallness of sigma}%
\end{equation}
we would have
\[
p\left(  u^{\ast}\right)  ^{p-1-\frac{qr}{s+1}}\leq1,
\]
hence $d_{1i}<0$, $i=1,2,3,\cdots$. So we have for any $d_{1}>0$, $A_{d_{1}
}=\emptyset$ and $N_{d_{1}}=0$.
\end{remark}

Finally, we give some sufficient conditions for the existence of \textit{a
priori} bounds uniform in $\tau$ of $\left(  \ref{Equation deformation}
\right)  $.

\begin{lemma}
\label{Lemma basic estimates for deformation equations}Let $\left(
u,v\right)  $ be a positive solution to $\left(  \ref{Equation deformation}
\right)  $. Then we have

\begin{enumerate}
\item[(i)]
\[
u\geq c\int_{\Omega}u,\quad v\geq c\int_{\Omega}v
\]
where $c>0$ is a constant independent of $\tau$.

\item[(ii)] For any $0<\gamma<\frac{n}{n-2}$,
\[
\int_{\Omega}u^{\gamma}\leq c\left(  \int_{\Omega}u\right)  ^{\gamma}
,\quad\int_{\Omega}v^{\gamma}\leq c\left(  \int_{\Omega}v\right)  ^{\gamma}
\]
where $c>0$ is a constant independent of $\tau$.

\item[(iii)]
\[
\tau\int_{\Omega}\frac{u^{p}}{v^{q}}\leq\int_{\Omega}u,\quad\tau\int_{\Omega
}\frac{u^{r}}{v^{s}}\leq\int_{\Omega}v,\quad\tau\int_{\Omega}\frac{u^{r}
}{v^{s+1}}\leq\left\vert \Omega\right\vert
\]
and if $\tau\in\left[  \frac{1}{2},1\right]  $,
\[
\int_{\Omega}v^{s+1}\leq2\int_{\Omega}u^{r}.
\]

\end{enumerate}
\end{lemma}

\begin{proof}
We refer the readers to the proofs of Lemmas
\ref{Lemma u,v larger than their average}, \ref{Lemma reverse Holder} and
\ref{Lemma trivial estimates}.
\end{proof}

When $\sigma>0$, \textit{a priori} lower bounds uniformly in $\tau$ can be
obtained using maximum principle.

\begin{proposition}
\label{Proposition Lower uniform bound when sigma>0}Let $\sigma>0$ and
$\left(  u,v\right)  $ be a positive solution to $\left(
\ref{Equation deformation}\right)  $, then
\[
u\geq c_{1},\quad v\geq c_{2}
\]
where $c_{1},c_{2}$ are positive constants depending only on $\sigma$ and
$\rho$.
\end{proposition}

\begin{proof}
Let $x^{\ast}\in\overline{\Omega}$ be a point such that
\[
u\left(  x^{\ast}\right)  =\inf_{x\in\overline{\Omega}}u\left(  x\right)  .
\]
Then we have at $x^{\ast}$, $\triangle u\geq0$, hence
\[
u\left(  x^{\ast}\right)  \geq\tau\left(  \frac{u^{p}}{v^{q}}+\sigma\right)
+\left(  1-\tau\right)  \rho\geq\tau\sigma+\left(  1-\tau\right)  \rho\geq
\min\left\{  \sigma,\rho\right\}  .
\]
Next, let $x^{\ast\ast}\in\overline{\Omega}$ be a point such that
\[
v\left(  x^{\ast\ast}\right)  =\inf_{x\in\overline{\Omega}}v\left(  x\right)
.
\]
Then we have at $x^{\ast\ast}$, $\triangle v\geq0$, hence
\[
v\left(  x^{\ast\ast}\right)  \geq\tau\frac{u^{r}}{v^{s}}+\left(  1-\chi
_{\tau}\right)  \rho.
\]
If $\tau\in\left[  0,\frac{1}{4}\right]  $, then we have
\[
v\left(  x^{\ast\ast}\right)  \geq\left(  1-\chi_{\tau}\right)  \rho\geq
\frac{1}{2}\rho.
\]
And if $\tau\in\left[  \frac{1}{4},1\right]  $, we have
\[
v\left(  x^{\ast\ast}\right)  \geq\tau\frac{u^{r}}{v^{s}}\geq\frac{1}{4}
\frac{u^{r}}{v^{s}},
\]
hence
\[
v\left(  x^{\ast\ast}\right)  \geq\left(  \frac{1}{4}u^{r}\right)  ^{\frac
{1}{s+1}}\geq\left(  \frac{1}{4}\left(  \min\left(  \sigma,\rho\right)
\right)  ^{r}\right)  ^{\frac{1}{s+1}}.
\]

\end{proof}

When $\sigma=0$, we will use the energy method in Section
\ref{Section energy 2} to obtain lower bounds.

\begin{proposition}
\label{Proposition delta}Assume $\sigma=0$, $r<\frac{n}{n-2}$ and there exists
$\delta\in\left(  0,1\right]  $ such that
\[
0<\frac{1-\delta}{r}+\frac{\delta}{p}<1,
\]
and
\[
\frac{\frac{\left(  1-\delta\right)  s}{r}+\frac{\delta q}{p}}{\frac
{r-1+\delta}{r}-\frac{\delta}{p}}<\frac{n}{n-2}\text{ or }\frac{\frac{\left(
1-\delta\right)  s}{r}+\frac{\delta q}{p}}{\frac{r-1+\delta}{r}-\frac{\delta
}{p}}\leq s+1.
\]
Then $\left(  \ref{Equation deformation}\right)  $ has a priori lower bounds
uniform in $\tau$.
\end{proposition}

\begin{proof}
First we assume $\tau\in\left[  \frac{1}{2},1\right]  $ and we will closely
follow the proof of Lemma \ref{Theorem a family of lower bounds}. Applying
H\"{o}lder inequality, we have
\begin{align*}
&  \int_{\Omega}u\leq\left(  \int_{\Omega}\frac{u^{r}}{v^{s}}\right)
^{\frac{1-\delta}{r}}\left(  \int_{\Omega}\frac{u^{p}}{v^{q}}\right)
^{\frac{\delta}{p}}\left(  \int_{\Omega}v^{\frac{\left(  \frac{1-\delta}
{r}\right)  s+\frac{\delta q}{p}}{\frac{r-1+\delta}{r}-\frac{\delta}{p}}
}\right)  ^{\frac{r-1+\delta}{r}-\frac{\delta}{p}}\\
\leq &  \left(  2\int_{\Omega}v\right)  ^{\frac{1-\delta}{r}}\left(
2\int_{\Omega}u\right)  ^{\frac{\delta}{p}}\left(  \int_{\Omega}
v^{\frac{\left(  \frac{1-\delta}{r}\right)  s+\frac{\delta q}{p}}
{\frac{r-1+\delta}{r}-\frac{\delta}{p}}}\right)  ^{\frac{r-1+\delta}{r}
-\frac{\delta}{p}}\\
\leq &  c\left(  \int_{\Omega}u\right)  ^{1+\frac{\delta}{p\left(  s+1\right)
}\left(  qr-\left(  p-1\right)  \left(  s+1\right)  \right)  }%
\end{align*}
where we have used part (ii) and part (iii) of Lemma
\ref{Lemma basic estimates for deformation equations}. Hence, we deduce
\[
\int_{\Omega}u\geq c.
\]
Applying part (i) of Lemma
\ref{Lemma basic estimates for deformation equations}, we have
\[
u\geq c
\]
where $c$ is a positive constant independent of $\tau\in\left[  \frac{1}
{2},1\right]  $. Applying maximum principle to the equation for $v$, we have
\[
\underline{v}\geq\left(  \tau\underline{u}^{r}\right)  ^{\frac{1}{s+1}}
\geq\left(  \frac{1}{2}c^{r}\right)  ^{\frac{1}{s+1}}.
\]
Next, we assume $\tau\in\left[  0,\frac{1}{2}\right]  $. Similar to the proof
of Proposition \ref{Proposition Lower uniform bound when sigma>0}, maximum
principle implies
\[
u\geq\frac{1}{2}\rho
\]
and
\[
v\geq\min\left\{  \frac{1}{2}\rho,\left(  \frac{1}{4}\left(  \frac{1}{2}
\rho\right)  ^{r}\right)  ^{\frac{1}{s+1}}\right\}  .
\]

\end{proof}

Once we have \textit{a priori} lower bounds, upper bounds can be obtained
using the method in Section \ref{Section maximum principle}.

\begin{proposition}
\label{Proposition upper uniform bound}Assume that we have positive lower
bounds for $\left(  \ref{Equation deformation}\right)  $ which is uniform in
$\tau$. Assume in addition that $\frac{p-1}{r}<1$. Then $\left(
\ref{Equation deformation}\right)  $ has a priori upper bounds uniform in
$\tau$.
\end{proposition}

\begin{proof}
Let $0<\lambda<\min\left\{  1,\frac{d_{2}}{2d_{1}}\right\}  $. At any point
$x^{\ast}\in\overline{\Omega}$ where
\[
\frac{u}{v^{\lambda}}\left(  x^{\ast}\right)  =\max_{x\in\overline{\Omega}
}\frac{u}{v^{\lambda}},
\]
following the proof of Lemma \ref{Lemma local maximum}, we have
\[
1-\tau\frac{u^{p-1}}{v^{q}}-\left(  \tau\sigma+\left(  1-\tau\right)
\rho\right)  u^{-1}-\frac{\lambda d_{1}}{d_{2}}\left(  1-\tau\frac{u^{r}
}{v^{s+1}}-\left(  1-\chi_{\tau}\right)  \rho v^{-1}\right)  \leq0.
\]
We rewrite the inequality into the form of
\begin{align*}
&  \left[  \left(  1-\frac{\lambda d_{1}}{d_{2}}-\frac{\tau}{2}\right)
-\left(  \tau\sigma+\left(  1-\tau\right)  \rho\right)  u^{-1}\right]
+\tau\left[  \frac{1}{2}+\frac{\lambda d_{1}}{d_{2}}\frac{u^{r}}{v^{s+1}
}-\frac{u^{p-1}}{v^{q}}\right] \\
\leq &  -\frac{\lambda d_{1}}{d_{2}}\left(  1-\chi_{\tau}\right)  \rho
v^{-1}\leq0.
\end{align*}
If
\[
\left(  1-\frac{\lambda d_{1}}{d_{2}}-\frac{\tau}{2}\right)  -\left(
\tau\sigma+\left(  1-\tau\right)  \rho\right)  u^{-1}\leq0,
\]
then we have
\[
u\left(  x^{\ast}\right)  \leq\frac{\tau\sigma+\left(  1-\tau\right)  \rho
}{1-\frac{\lambda d_{1}}{d_{2}}-\frac{\tau}{2}}\leq\frac{\sigma+\rho}{\frac
{1}{2}-\frac{\lambda d_{1}}{d_{2}}},
\]
hence, using the positive lower bound for $v$,
\[
\frac{u}{v^{\lambda}}\left(  x^{\ast}\right)  \leq\frac{u\left(  x^{\ast
}\right)  }{\overline{v}^{\lambda}}\leq c
\]
where $c$ is a positive constant independent of $\tau$.

Otherwise, we have
\[
\tau\left[  \frac{1}{2}+\frac{\lambda d_{1}}{d_{2}}\frac{u^{r}}{v^{s+1}}
-\frac{u^{p-1}}{v^{q}}\right]  <0.
\]
Since $\frac{p-1}{r}<1$ and $\left(  \ref{Equation pqrs}\right)  $ holds, we
can choose $\varepsilon>0$ such that
\[
\frac{p-1}{r}<\frac{q-\varepsilon}{s+1}<1,
\]
and we further assume
\[
\lambda<\frac{s+1-\left(  q-\varepsilon\right)  }{r-\left(  p-1\right)  }.
\]
Let
\[
a_{\varepsilon}=\frac{\left(  q-\varepsilon\right)  -\lambda\left(
p-1\right)  }{s+1-\lambda r},
\]
then it is easy to verify $a_{\varepsilon}\in\left(  0,1\right)  $. Using
Young's inequality,
\begin{align*}
&  \frac{1}{2}+\frac{\lambda d_{1}}{d_{2}}\frac{u^{r}}{v^{s+1}}<\frac{u^{p-1}
}{v^{q}}\leq\underline{v}^{-\varepsilon}\frac{u^{p-1}}{v^{q-\varepsilon}}\\
=  &  \underline{v}^{-\varepsilon}\left(  \frac{u^{r}}{v^{s+1}}\right)
^{a_{\varepsilon}}\left(  \left(  \frac{u}{v^{\lambda}}\right)  ^{-\frac
{\left(  q-\varepsilon\right)  r-\left(  p-1\right)  \left(  s+1\right)
}{s+1-\lambda r-\left(  q-\varepsilon-\lambda\left(  p-1\right)  \right)  }
}\right)  ^{1-a_{\varepsilon}}\\
\leq &  \frac{\lambda d_{1}}{d_{2}}\frac{u^{r}}{v^{s+1}}+c\left(  \frac
{u}{v^{\lambda}}\right)  ^{-\frac{\left(  q-\varepsilon\right)  r-\left(
p-1\right)  \left(  s+1\right)  }{s+1-\lambda r-\left(  q-\varepsilon
-\lambda\left(  p-1\right)  \right)  }},
\end{align*}
hence we again have
\[
\frac{u}{v^{\lambda}}\left(  x^{\ast}\right)  \leq c
\]
where $c$ is a positive constant independent of $\tau$. Finally, let
$x^{\ast\ast}\in\overline{\Omega}$ be such that
\[
u\left(  x^{\ast\ast}\right)  =\max_{x\in\overline{\Omega}}u.
\]
From maximum principle, we have at $x^{\ast\ast}$,
\begin{align*}
&  u\leq\tau\left(  \frac{u^{p}}{v^{q}}+\sigma\right)  +\left(  1-\tau\right)
\rho\leq\frac{u^{p}}{v^{q}}+\sigma+\rho\\
\leq &  c^{p}v^{p\lambda-q}+\sigma+\rho.
\end{align*}
If we choose $\lambda$ so that $\lambda<\frac{q}{p}$, we have
\[
\max_{x\in\overline{\Omega}}u=u\left(  x^{\ast\ast}\right)  \leq
c^{p}\underline{v}^{p\lambda-q}+\sigma+\rho
\]
which is a bound independent of $\tau$. Finally, at the point where $v$
achieves its maximum, we have
\[
v\leq\tau\frac{u^{r}}{v^{s}}+\left(  1-\chi_{\tau}\right)  \rho\leq
\frac{\overline{u}^{r}}{\underline{v}^{s}}+\rho.
\]

\end{proof}

Combining \textit{a priori} estimates in Propositions
\ref{Proposition Lower uniform bound when sigma>0}, \ref{Proposition delta}
and \ref{Proposition upper uniform bound} with Theorem \ref{Theorem existence}
, we have

\begin{theorem}
\label{Theorem existence detailed}Assume that $\sigma>0$ or the assumptions in
Proposition \ref{Proposition delta} hold. Assume in addition that $\frac
{p-1}{r}<1$. If $d_{1}\neq d_{1i}$, $i=1,2,3,\cdots$, and $N_{d_{1}}$ is odd,
then there exists at least one nontrivial solution to $\left(
\ref{Equation stationary}\right)  $.
\end{theorem}

\section*{Acknowledgements}

We wish to thank the anonymous referee for pointing out a flaw in the original
proof of Lemma \ref{Lemma index 1} and providing several additional
references. The research is supported in part by the NSF.

\end{document}